\documentclass[journal]{IEEEtran}
\usepackage{graphicx}
\usepackage{epstopdf}
\usepackage{epsfig}
\usepackage{amsfonts}
\usepackage{amstext}
\usepackage{color}
\usepackage{empheq}
\usepackage{graphicx,colordvi}
\usepackage{amssymb,amsmath,amsopn,amsfonts,graphicx}
\usepackage{CJK}
\usepackage{float}
\usepackage[justification=centering]{caption}

\newtheorem{condition**}{A*}
\newtheorem{condition***}{C*}
\newtheorem{condition*}{C}

\newtheorem{proposition}{Proposition}[section]

\newtheorem{lemma}{Lemma}[section]
\newtheorem{theorem}{Theorem}[section]
\newtheorem{remark}{Remark}[section]

\date{}

\begin{document}

\title{ Social Optima in Robust Mean Field LQG Control: From Finite to Infinite Horizon}
\author{Bing-Chang Wang,~\IEEEmembership{Member, IEEE }
        Jianhui Huang,~\IEEEmembership{ }
        Ji-Feng~Zhang,~\IEEEmembership{Fellow, IEEE}
\thanks{The work was partially supported by the National Natural Science Foundation of China under Grants 61773241, 61877057, and partially supported by the PolyU G-YL04, RGC Grants 502412, 15300514. This paper was presented in a preliminary form at the 2017 Asian Control Conference.}

\thanks{B. C. Wang is with School of Control Science and Engineering, Shandong University, Jinan 250061, P. R. China (e-mail: bcwang@sdu.edu.cn).}

\thanks{J. Huang is with the Department of Applied Mathematics,
The Hong Kong Polytechnic University, Hong Kong (e-mail: majhuang@polyu.edu.hk).}

\thanks{J. F. Zhang is with the Key Laboratory of Systems and Control, Academy of Mathematics and Systems Science, Chinese Academy of Sciences, Beijing 100190, China; and and also with School of Mathematical Sciences, University of Chinese Academy of
Sciences, Beijing 100149, China (e-mail: jif@iss.ac.cn). }}
\maketitle

{\bf Abstract} -- This paper studies social optimal control of mean field LQG (linear-quadratic-Gaussian) models with uncertainty. Specially, the uncertainty is represented by a uncertain drift which is common for all agents. A robust optimization approach is applied by assuming all agents treat the uncertain drift as an adversarial player. In our model, both dynamics and costs of agents are coupled by mean field terms, and both finite- and infinite-time horizon cases are considered. By examining social functional variation and exploiting person-by-person optimality principle, we construct an auxiliary control problem for the generic agent via a class of forward-backward stochastic differential equation system. By solving the auxiliary problem and constructing consistent mean field approximation, a set of decentralized control strategies is designed and shown to be asymptotically optimal.


\section{Introduction}

\subsection{Background and Motivation}

Mean field games and control have drawn increasing attention in many fields, including system control, applied mathematics and economics \cite{BFY13, C14, GS13}. The mean field game involves a very large population of small interacting players with the feature that while the influence of each one is negligible, the impact of the overall population is significant. 
By now, mean field games and control have been intensively studied in the linear-quadratic-Gaussian (LQG) framework \cite{HCM03, HCM07, LZ08, WZ13, MB17}, and
there is a large body of works on nonlinear models \cite{HMC06, LL07, CD13}. Huang \emph{et al.} designed $\epsilon$-Nash equilibrium strategies for LQG mean field games with discount costs based on the proposed Nash certainty equivalence (NCE) approach \cite{HCM03, HCM07}. The NCE approach was then applied to the cases with long run average costs \cite{LZ08} and with Markov jump parameters \cite{WZ12}, respectively. Lasry and Lions independently introduced the model of mean field games and studied the well-posedness problem of the limiting partial differential equations \cite{LL07}. For further literature, readers are referred to \cite{H10, WZ12, WZ14} on mean field games with major players in continuous- or discrete-time,
\cite{CD13} on probabilistic analysis of mean field games, and
\cite{weintraub2008markov} on the oblivious equilibrium in dynamic games.

Besides noncooperative games, social optima in mean field models have also drawn much attention. The social optimum control refers to that all the players cooperate to optimize the common social cost---the sum of individual costs, which is usually regarded as a type of team decision problem \cite{H80}.  Huang \emph{et al.} considered social optima in mean field LQG control, and provided an asymptotic team-optimal solution \cite{HCM12}. Wang and Zhang investigated a mean field social optimal problem where a Markov jump parameter appears as a common source of randomness \cite{WZ17}.
Also, see \cite{HN16} for social optima in mixed games, \cite{AM15} for team-optimal control with finite population and partial information, and \cite{LZZ16} for social optima of static mean field games.

Mathematical models can only be approximations of the real world. Actually, some parts of a model may be inexact. Thus, it is worthwhile to study the mean field control with model uncertainty \cite{BB95}. The works \cite{HH13, HH16, WH17} investigated the mean field games and control with a global uncertainty term.  The ``hard constraint" case (the disturbance is specified with a bound) was considered in \cite{HH13} under which the substantial difficulty arises after the Lagrange multiplier is introduced. Authors in \cite{HH16, WH17} adopted the ``soft constraint" approach (\cite{BB95,BES03,E06}) by removing the bound of the disturbance while the effort is penalized in the cost function. The works \cite{TBT13, MB17} considered the case that each agent is paired with the local disturbance as an adversarial player, and provided an $\epsilon$-Nash equilibrium by tackling a Hamilton-Jacobi-Isaacs equation combined with a fixed-point analysis.

\subsection{Challenge and Contribution}

This paper investigates mean field LQG social optimum control with a common uncertain drift, where both dynamics and costs of agents involve mean field coupled terms. To address the model uncertainty, a minus quadratic penalty term of drift is incorporated into the cost functional.
There exist some substantial challenges in studying the problem.
First, the socially optimal control with respect to drift uncertainty is a high-dimensional optimization problem with indefinite state weights.
The corresponding convexity condition is very hard to verify. Second, by social variational derivation, the resulting limit system is governed by a controlled  forward-backward stochastic differential equation (FBSDE). To design
decentralized strategies, we need to solve the auxiliary optimal control problem subject to an FBSDE system.
Meanwhile, the asymptotic optimality analysis is 
different
from the general mean field LQG problems.
Third, for the social optimum problem in the infinite horizon, we are faced with tackle infinite-horizon FBSDEs and the relevant optimal control problems.

In this paper, the social optimum control for the robust mean field LQG model is tackled by using stochastic maximum principle \cite{YZ99, ZQ16, ZX18}. For the finite-horizon problem, we first obtain some low-dimensional convexity conditions and a set of FBSDEs by analyzing the variation of the centralized maximization cost to drift uncertainty. With the help of the Riccati equation, we further obtain a feedback type of the ``worst-case" drift for the social optimum problem. Next, we construct an auxiliary optimal control problem based on the social variational derivation and the person-by-person optimality principle. By solving the auxiliary problem combined with consistent mean field approximations, a set of decentralized control laws is designed and further shown to be asymptotically robust social optimal by perturbation analysis. Finally, from asymptotic analysis to FBSDEs we design decentralized strategies and show their robust optimality for the infinite-horizon social optimum problem.

The main contributions of the paper are summarized as follows.
(i)
Social optimum control is studied for mean field models with a common uncertain drift, where coupled terms are included in both costs and dynamics of agents. By FBSDE and Riccati equation approaches we design a set of decentralized feedback control laws.
(ii) By examining the social cost variation, we give low-dimensional convexity conditions and asymptotic convexity analysis for robust social optimum problems.
(iii) From consistency requirements in mean field approximations, a system of differential equations is derived. The existence condition of solutions to consistency equations is characterized by a Riccati equation, instead of a fixed-point analysis.
(iv) From the perturbation analysis to FBSDE, the decentralized strategies are shown to have asymptotic robust optimality.
(v) By analyzing the asymptotic behavior of FBSDE, decentralized strategies for the infinite-horizon problem are designed and further shown be robust social optimal.

\subsection{Organization and Notation}
The organization of the paper is as follows. In Section II, we consider the finite-horizon social optimization problem with drift uncertainty. By variational analysis, the centralized control with respect to drift uncertainty is obtained. Then an auxiliary optimal control problem is constructed based on person-by-person optimality. By solving this problem combined with consistent mean field approximations, a set of decentralized strategies is designed and further proved to be robust social optimal. Section III tackles the infinite-horizon social optimum problem. In Section IV, a numerical example is provided to verify the result.
Section V concludes the paper.

{\it Notation:} Suppose that $(\Omega,
\mathcal F, \{\mathcal F_t\}_{0\leq t\leq T}, \mathbb{P})$ is a complete
filtered probability space.
  Throughout this paper, we denote by 
   $\otimes$ the Kronecker product, $I_m$ m-dimensional identity matrix ($I_n$ abbreviated as $I$).
 We use  $\|\cdot\|$ to denote the norm of a  Euclidean space, or the Frobenius norm for matrices. For a symmetric matrix $Q$ and a vector $z$, $\|z\|_Q^2= z^TQz$; for two vectors $x,y$, $\langle x,y\rangle=x^Ty$. For a matrix (vector) $M$, $M^T$ denotes its transpose, $M>0$ means that $M$ is positive definite.
 Let $L^{2}_{\mathcal{F}}(0, T; \mathbb{R}^{k})$ denote the space of all $\mathbb{R}^{k}$-valued $\mathcal{F}_t$-progressively measurable processes $x(\cdot)$ satisfying $\mathbb{E}\int_0^{T}\|x(t)\|^{2}dt<\infty$, and $L^{2}_{\mathcal{F},\frac{\rho}{2}}(0, \infty; \mathbb{R}^{k})$ denote the space of all $\mathbb{R}^{k}$-valued $\mathcal{F}_t$-progressively measurable processes $x(\cdot)$ satisfying  $\mathbb{E}\int_0^{\infty}e^{-\rho t}\|x(t)\|^2dt<\infty .$
 $C([0,T],\mathbb{R}^k)$ is the space of all $\mathbb{R}^k$-valued functions defined on $[0,T]$ which are continuous; $C_{\rho/2}([0,\infty), \mathbb{R}^k)$ is a subspace of $C([0,\infty),\mathbb{R}^k)$ which is given by $\{f|\int_0^{\infty}e^{-\rho t}\|f(t)\|^2dt<\infty  \}.$
 For convenience of  presentation, we use $C$ (or $ C_1, C_2, \ldots$) to denote a generic constant which 
 may vary from place to place.


\section{Mean Field Social Control over a Finite Horizon}

Consider a large population systems with $N$ agents. The $i$th agent evolves by the following stochastic differential equation:
\begin{align}\label{eq1}
dx_i(t)=& [Ax_i(t)+Bu_i(t)+Gx^{(N)}(t)+f(t)]dt\cr
&+\sigma dW_i(t),\  1\leq i\leq N,
\end{align}
where $x_i\in
\mathbb{R}^n$ and $u_i\in\mathbb{R}^r$ are the state and the input of agent $i$, respectively. $x^{(N)}(t)=\frac{1}{N}\sum_{j=1}^Nx_j(t)$. $\{W_i, 1\leq i\leq N\}$ are a sequence of mutually independent $d$-dimensional Brownian motions. $f\in  L^{2}_{\mathcal{F}}(0, T; \mathbb{R}^{n})$ is an unknown disturbance, which reflects the effect imposed to each agent by the eternal environment.
 The cost function of agent $i$ is given by
\begin{align}\label{eq2}
J_i^{\rm F}(u)=&\frac{1}{2}\mathbb{E}\int_0^{T}\Big\{
 \big\|x_i(t)-\Gamma x^{(N)}(t)-\eta\big\|^2_{Q}\cr
&+\|u_i(t)\|^2_{R_1}-\|f(t)\|^2_{R_2}\Big\}dt\cr
&+\frac{1}{2}\mathbb{E}\|x_i(T)\|^2_H,
\end{align}
where $Q, R_1,R_2, H\in \mathbb{R}^{n\times n}$ are symmetric, $\Gamma\in \mathbb{R}^{n\times n}$ and $\eta\in \mathbb{R}^n$.
  $u=\{u_1,
\ldots, u_N\}$. Take $\{{\mathcal F}_t\}_{0\leq t\leq T}$ as the natural filtration generated by the $Nd$-dimensional Brownian motion $(W_1,\cdots, W_N)$.
The decentralized control set is given by
\begin{align}
{\cal U}_{i}^{\rm F}=\Big\{u_i\ |\
  & u_i(t)\in \sigma(x_i(s), 0\leq s\leq t), \cr
  & \mathbb{E}\int_0^T\|u_i(t)\|^2dt<\infty \Big\}.\nonumber
\end{align}
For comparison, define the centralized control set as
\begin{align*}
{\cal U}_{c}^{\rm F} =\Big\{u_i|&\
  u_i(t)\in \sigma(x_i(0),W_i(s), 0\leq s\leq t,\ 1\leq i\leq N),\cr
   & \mathbb{E}\int_0^T\|u_i(t)\|^2dt<\infty \Big\}.
\end{align*}
Denote $J_{\rm soc}^{\rm F}(u)=\sum_{i=1}^NJ_i^{\rm F}(u)$.
Let the social cost under the worst-case disturbance be
$$J_{\rm soc}^{\rm wo}({u})=
\sup_{f\in {\mathcal U}_c^{\rm F}}J_{\rm soc}^{\rm F}({u},f).$$
\textbf{Problem (PF):} Seek a set of decentralized control laws $(\hat{u}_1,\cdots,\hat{u}_N)$ to minimize the social cost under the worst-case disturbance
for System (\ref{eq1})-(\ref{eq2}), i.e.,
$\inf_{u_i\in {\cal U}_{i}^{\rm F}}J_{\rm soc}^{\rm wo}({u})$.

\begin{remark}
 Different from \cite{HH16, TBT13}, we assume the disturbance $f$ is a common stochastic process. $f$ may stand for the impact from tax, subsidy or natural disaster. In this case, agents may be pessimistic to suppose the disturbance would use the information of all the agents to play against them.
 \end{remark}

We make the following assumptions.

\textbf{(A0)} $\{x_i(0)\}$ are independent random variables with the same mathematical expectation. $x_i(0)=x_{i0}$, $\mathbb{E}x_i(0)=\bar{x}_0,
1\leq i \leq N$. There exists a constant $C_0$ such that $\max_{1\leq i \leq N}\mathbb{E}\|x_{i0}\|^2<C_0$.
Furthermore, $\{x_{i0}, i=1,...,N\}$ and
$\{W_i, i=1,...,N\}$ are independent of each other.

\textbf{(A1)} $Q\geq0$, $R_1>0$, $R_2>0$, and $H\geq0$.

From now on, the time variable $t$ might be suppressed if necessary and no confusion occurs.

\subsection{The Control Problem with Respect to Model Uncertainty}

Let $u_i=\check{u}_i\in {\cal U}_{c}^{\rm F},  i=1,\cdots, N$ be fixed.
The optimal control problem with respect to drift uncertainty is as follow:
$$\hbox{\textbf{(P1)} maximize}_{f\in {\mathcal U}_{c}^{\rm F}}J_{\rm soc}^{\rm F}(\check{u},f).$$
Clearly, (P1) is equivalent to the following problem:
$$\hbox{\textbf{(P1$^\prime$)} minimize}_{f\in{\mathcal U}_{c}^{\rm F}}\check{J}_{\rm soc}^{\rm F}(f),$$
where
\begin{align*}
  \check{J}_{\rm soc}^{\rm F}(f)=&\frac{1}{2}\sum_{i=1}^N
\mathbb{E}\int_0^{T}\Big\{-\big\|x_i
   -\Gamma x^{(N)}-\eta\big\|^2_{Q}
+\|f\|^2_{R_2}\Big\}dt\cr
&-\frac{1}{2}\mathbb{E}\|x_i(T)\|^2_H.
\end{align*}

Let  $\textbf{x}=(x_1^T,\cdots,x_N^T)^T$, $\textbf{u}=(u_1^T,\cdots,u_N^T)^T$, $\textbf{1}=(1,\cdots,1)^T$, $\textbf{W}=(W_1^T,\cdots,W_N^T)^T$,
$\textbf{A}=diag(A, \cdots,A)$, $\textbf{B}=diag(B, \cdots,B)$,
${\hat \sigma}=diag(\sigma, \cdots,\sigma)$,  $\textbf{H}=diag(H, \cdots,H)$, $\hat{\textbf{Q}}=diag\{Q,\cdots,Q\}-\frac{1}{N}\textbf{11}^T\otimes \Psi,$ and $\hat{\bf{\eta}}=\textbf{1}\otimes \bar{\eta}$, where $\Psi\stackrel{\Delta}{=}\Gamma^TQ+Q\Gamma-\Gamma^TQ\Gamma$  and $\bar{\eta}\stackrel{\Delta}{=}Q\eta-\Gamma^T Q\eta$.
We can write Problem (P1$^{\prime}$) as to
minimize
\begin{align*}
\check{J}_{\rm soc}^{\rm F}(f)=&\frac{1}{2}\mathbb{E}\int_0^{T}\Big(-\textbf{x}^T\hat{\textbf{Q}}\textbf{x}+2\hat{\bf{\eta}}^T\textbf{x}
+{N}f^TR_2f\Big)dt\cr
&-\frac{1}{2}\mathbb{E}[\textbf{x}^T(T)\textbf{Hx}(T)],
\end{align*}
subject to
$$
d\textbf{x}(t)=\check{\textbf{A}}\textbf{x}(t)dt+\textbf{Bu}(t)dt+ \textbf{1}\otimes f(t)dt+{\hat{\sigma}} d\textbf{W}(t).
$$
where $\check{\textbf{A}}\stackrel{\Delta}{=}{\textbf{A}}+\frac{1}{N}(\textbf{11}^T\otimes G).$

For the further existence analysis, we introduce the following assumptions:

\textbf{(A2)} Problem (P1$^{\prime}$) is convex in $f$;

\textbf{(A2$^{\prime}$)} Problem (P1$^{\prime}$) is uniformly convex in $f$.

Below are some necessary and sufficient  conditions to ensure (A2) or (A2$^{\prime}$). 

\begin{proposition}\label{prop1}
The following statements are equivalent:

(i) Problem (P1$^{\prime}$) is convex in $f$.

(ii)  For any $f\in {\mathcal U}_{c}^{\rm F}$,
 $$\int_0^{T}\Big(-\textbf{y}^T\hat{\textbf{Q}}\textbf{y}
+{N}f^TR_2f\Big)dt-\|\textbf{y}(T)\|^2_{\textbf{H}}\geq 0,$$
where $\textbf{y}\in \mathbb{R}^{nN}$ satisfies
\begin{equation*}\label{eq3a}
   d\textbf{y} = (\check{\textbf{A}}\textbf{y}+\textbf{1}\otimes f)dt,\
   \textbf{y}(0)= 0.
\end{equation*}

(iii) For any $f\in {\mathcal U}_{c}^{\rm F}$,
$$\mathbb{E}\int_0^T\Big\{-\big\|(I-\Gamma)y_i\big\|^2_{Q}+\|f\|^2_{R_2}\Big\}dt-\|y_i(T)\|_H^2\geq 0, $$
where for $i=1$, 2, $\cdots$, $N$, $y_i$ satisfies
\begin{equation}\label{eq4}
dy_i=[A y_i+Gy^{(N)}+f]dt,\quad y_i(0)=0.
\end{equation}
\end{proposition}

\emph{Proof.} (i) $\Leftrightarrow$ (ii) is given in \cite{HH16, LZ99}.
From (\ref{eq4}), we have $y_1=y_2=\cdots=y_N=y^{(N)}$.
Thus,
\begin{align}\label{eq4b}
  &\int_0^{T}\Big(-\textbf{y}^T\hat{\textbf{Q}}\textbf{y}
+{N}f^TR_2f\Big)dt-\|\textbf{y}(T)\|^2_{\textbf{H}}\cr
=&\sum_{i=1}^N\!\!\int_0^T\!\!\Big(\!\!-\big\| y_i-\Gamma y_i\big\|^2_{Q}+\|f\|^2_{R_2}\Big)dt\!-\!\!\sum_{i=1}^N\!\|y_i(T)\|_H^2\cr
=&N\!\Big[\!\int_0^T\!\!\Big(\!\!-\big\|(I\!-\!\Gamma)y_i\big\|^2_{Q}\!+\!\|f\|^2_{R_2}\Big)dt\!-\!\|y_i(T)\|_H^2\Big],\cr
\end{align}
which implies that (ii) is equivalent to (iii). $\hfill \Box$

\begin{proposition}\label{prop1b}
The following statements are equivalent:

(i) Problem (P1$^{\prime}$) is uniformly convex in $f$.

(ii) There exists $\delta>0$
such that
 $$\int_0^{T}\!\!\Big(\!-\textbf{y}^T\hat{\textbf{Q}}\textbf{y}
+{N}f^TR_2f\Big)dt-\|\textbf{y}(T)\|^2_{\textbf{H}}\geq \delta\!\int_0^T\!\!\|f\|^2dt.$$

(iii) The equation
$$
\dot{\textbf{P}}+\check{\textbf{A}}^T\textbf{P}+\textbf{P}\check{\textbf{A}}-\hat{\textbf{Q}}-
\textbf{P}(\textbf{1}\otimes I) (NR_2)^{-1}(\textbf{1}^T\otimes I)\textbf{P}=0,
$$
with $\textbf{P}(T)=-\textbf{H}$ admits a solution in $C([0,T];\mathbb{R}^{nN})$.

(iv) The following equation admits a solution in $C([0,T];\mathbb{R}^{n})$,
\begin{align*}
&\dot{P}+(A+G)^TP+P(A+G)-PR_2^{-1}P\cr
&-(I-\Gamma)^TQ(I-\Gamma)=0,\quad P(T)=-H.
\end{align*}

(v) For any $t \in [0, T],$
$
\det[(0, I)e^{\mathbb{A}t}(0, I)^{T}]>0,
$
where
$$\mathbb{A}=\left(\begin{array}{cc} A+G+R^{-1}_2 H & -R^{-1}_2 \\[0.3cm]
\mathbb{A}_{21} & -(A+G+R^{-1}_2 H)^{T}
\end{array}\right)
$$
with
$
\mathbb{A}_{21}=HR^{-1}_2H\!+\!(I\!-\!\Gamma)^TQ(I\!-\!\Gamma)\!+\!(A\!+\!G)^{T}H\!+\!H(A\!+\!G).
$
\end{proposition}

\emph{Proof.} (i)$\Leftrightarrow$(ii) is implied from \cite{HH16, LZ99}. (i)$\Leftrightarrow$(iii) is given by Theorem 4.5 of \cite{SLY16}. \textcolor{blue}{By (\ref{eq4b}) and (ii), we have
$$\mathbb{E}\int_0^T\!\!\!\big\{-\|y_i\|^2_{Q(I-\Gamma)}+\|f\|^2_{R_2}\big\}dt-\|y_i(T)\|_H^2\geq \frac{\delta}{N}\mathbb{E}\int_0^T\!\!\|f\|^2dt. $$
By \cite[Theorem 4.5]{SLY16}, we obtain (ii)$\Leftrightarrow$(iv), which further implies (i)$\Leftrightarrow$(iv).}
(iv)$\Leftrightarrow$(v) is given by \cite{HH16, MY99}.
$\hfill \Box$

By examining the variation of $ \check{J}_{\rm soc}^{\rm F}$, we first
obtain the necessary and sufficient conditions for
the existence of centralized optimal control of (P1).

\begin{theorem}\label{thm1}
 (P1$^{\prime}$) has a minimizer in ${\mathcal U}_{c}^{\rm F}$ if and only if (A2) holds and the following equation system admits a set of 
 adapted solutions $(x_i,p_i, i=1,\cdots,N)$: 
\begin{equation}\label{eq4a}
\left\{
\begin{aligned}
 d\check{x}_i=&(A\check{x}_i+B\check{u}_i+{G}\check{x}^{(N)}\!-\!{R_2^{-1}}
 \check{p}_j^{(N)})dt\!+\!\sigma dW_i,\\
d\check{p}_i=&-[A^T\check{p}_i+{G^T}\check{p}_j^{(N)}-Q\check{x}_i+{\Psi }\check{x}^{(N)}+\bar{\eta} ]dt\cr
&+\sum_{j=1}^N\beta_i^jdW_j,\\
\check{x}_i(0)&={x_{i0}},\ \check{p}_i(T)=-H\check{x}_i(T),\ 1\leq i\leq N,
\end{aligned}\right.
\end{equation}
where $p^{(N)}=\frac{1}{N}\sum_{j=1}^Np_j$, 
and furthermore the minimizer is $\check{f}=-{R_2^{-1}}p^{(N)}$.
\end{theorem}

{\it Proof.}
Suppose that $\check{f}=-R_2^{-1}p^{(N)}$, where $p_i, i=1,\cdots,N$ are a set of solutions to the equation system
\begin{equation}\label{eq6}
dp_i=\alpha_idt+\beta_i^idW_i+\sum_{j\not =i}\beta_i^jdW_j, \ p_i(T)=-Hx_i(T),
\end{equation}
where $i=1$, $\cdots$, $N$; $\alpha_i$ and $\beta_i^j$ are to be determined. 
 Denote by $\check{x}_i$ the state of agent $i$ under the control $\check{u}_i$ and the drift $\check{f}$. For any $f\in {\mathcal U}_{c}^{\rm F} $ and $\varepsilon\in \mathbb{R}$, let $f^{\varepsilon}=\check{f}+\varepsilon f$. Let $ x_i^{\varepsilon}$ be the solution of the following perturbed state equation
$$ dx_i^{\varepsilon}=\big(Ax_i^{\varepsilon}+B\check{u}_i+\check{f}+\varepsilon f+\frac{G}{N}\sum_{i=1}^Nx^{\varepsilon}_i\big)dt+\sigma dW_i,
$$
with $x_i^{\varepsilon}(0)=x_{i0}$, $i=1$, 2, $\cdots$, $N$.

Let $y_i=(x_i^{\varepsilon}-\check{x}_i)/\varepsilon$, and $y^{(N)}=\sum_{i=1}^Ny_i/{N}$.
It can be verified that
 ${y}_i$ satisfies (\ref{eq4}).
Then, by It\^{o}'s formula,
\begin{align}\label{eq7}
 &-\mathbb{E}\langle H \check{x}_i(T),y_i(T) \rangle\cr
=&\mathbb{E}[\langle p_i(T),y_i(T)\rangle-\langle p_i(0),y_i(0)\rangle]\cr
=&\mathbb{E}\int_0^T \big[\langle \alpha_i,y_i\rangle+\langle p_i,Ay_i+Gy^{(N)}+f\rangle\big] dt.
\end{align}
We have
\begin{equation}\label{eq5a}
\begin{aligned}
 &\check{J}_{\rm soc}^{\rm F}(\check{f}+\varepsilon f)-\check{J}_{\rm soc}^{\rm F}(\check{f})=\varepsilon \Lambda_1+\frac{\varepsilon^2}{2}\Lambda_2
\end{aligned}
\end{equation}
where
\begin{align*}
\Lambda_1\stackrel{\Delta}{=}&\sum_{i=1}^N\mathbb{E}\int_0^T\!\! \Big[\big\langle\!\!-\!Q\big(\check{x}_i\!-\!(\Gamma\check{x}^{(N)}\!+\!\eta)\big),y_i\!-\!\Gamma y^{(N)}\big\rangle\cr
  & +
\langle R_2 \check{f},f\rangle \Big]dt-\sum_{i=1}^N\mathbb{E}\langle H \check{x}_i(T),y_i(T) \rangle,\cr
\Lambda_2\stackrel{\Delta}{=}&\sum_{i=1}^N\mathbb{E}\int_0^T\Big\{-\big\|y_i
   -\Gamma y^{(N)}\big\|^2_{Q}
+\|f\|^2_{R_2}\Big\}dt.
\end{align*}
Note that
\begin{align*}
&\sum_{i=1}^N\mathbb{E}\int_0^T\big\langle-Q\big(\check{x}_i-(\Gamma\check{x}^{(N)}+\eta)\big),\Gamma y^{(N)}\big\rangle dt\cr
=&\mathbb{ E}\int_0^T\!\!\big\langle\!-\!\Gamma^TQ  \sum_{i=1}^N\big(\check{x}_i-(\Gamma\check{x}^{(N)}+\eta)\big),\frac{1}{N}  \sum_{j=1}^Ny_j\big\rangle  dt\cr
 =& \sum_{j=1}^N \mathbb{E}\int_0^T\!\!\big\langle- \frac{\Gamma^TQ}{N} \sum_{i=1}^N\big(\check{x}_i-(\Gamma\check{x}^{(N)}+\eta)\big), y_j\big\rangle  dt\cr
 =& \sum_{j=1}^N \mathbb{E}\int_0^T\!\!\big\langle- {\Gamma^TQ} \big((I-\Gamma)\check{x}^{(N)}-\eta\big), y_j\big\rangle  dt.
\end{align*}
From (\ref{eq7}), one can obtain that
\begin{align*}
  \Lambda_1=&\mathbb{E}\sum_{i=1}^N\int_0^T\Big[\big\langle-Q\big(\check{x}_i-(\Gamma\check{x}^{(N)}+\eta)\big),
  \cr
  &y_i-\Gamma y^{(N)}\big\rangle +
  \langle R_2\check{f},f\rangle \Big]dt\cr
  +&\sum_{i=1}^N \mathbb{E}\int_0^T\big[\langle \alpha_i,y_i\rangle+\langle p_i,Ay_i+Gy^{(N)}+f\rangle\big] dt\cr
  =&\mathbb{E}\int_0^T\Big\langle NR_2\check{f}+\sum_{i=1}^Np_i,f\Big\rangle dt\cr
  +&\sum_{i=1}^N\mathbb{E}\int_0^T\Big\langle-Q\big(\check{x}_i-(\Gamma\check{x}^{(N)}+\eta)\big)
  +A^Tp_i\cr
  +& {\Gamma^TQ} \big((I-\Gamma)\check{x}^{(N)}-\eta\big)+\alpha_i+{G^T}p^{(N)}, y_i\Big\rangle dt,
\end{align*}
From (\ref{eq5a}), $\check{f}$ is a minimizer to Problem (P1$^{\prime}$) if and only if
$\Lambda_2\geq0$ and $\Lambda_1=0 $.
By Proposition \ref{prop1}, $\Lambda_2\geq0$ if and only if (A2) holds. $\Lambda_1=0 $ is equivalent to
\begin{align*}
  \alpha_i=&-\big[A^Tp_i+{G^T}p^{(N)}+\Gamma^TQ \big((I-\Gamma)\check{x}^{(N)}-\eta\big)\cr
  &-Q\big(\check{x}_i-(\Gamma\check{x}^{(N)}+\eta))\big],\cr
\check{f}=&-{R_2^{-1}}p^{(N)}.
\end{align*}
Thus, we have the optimality system (\ref{eq4a}).
Namely, $\Lambda_1=0 $ if and only if (\ref{eq4a}) admits a solution $(\check{x}_i,\check{p}_i, i=1,\cdots,N)$. \hfill$\Box$


Let 
$\check{u}^{(N)}=\frac{1}{N}\sum_{i=1}^N\check{u}_i$,  and $\check{p}^{(N)}=\frac{1}{N}\sum_{i=1}^N\check{p}_i$. It follows from (\ref{eq4a}) that
\begin{equation}\label{eq10}
\left\{
\begin{aligned}
 d\check{x}^{(N)}=&\big((A+G)\check{x}^{(N)}+B\check{u}^{(N)}-R_2^{-1}\check{p}^{(N)}\big)dt
 \\&+\frac{1}{N}\sum_{i=1}^N\sigma dW_i,\\
d\check{p}^{(N)}=&-\!\Big[(A+G )^T\check{p}^{(N)}+(\Psi-Q)\check{x}^{(N)}
+\bar{\eta}\Big]dt\\&
+\frac{1}{N}\sum_{i=1}^N\sum_{j=1}^N{\beta}_i^jdW_j,\\
\check{x}^{(N)}(0)&=\frac{1}{N}\sum_{i=1}^Nx_{i0},\ \check{p}^{(N)}(T)=-H \check{x}^{(N)}(T).
\end{aligned}\right.
\end{equation}
\begin{proposition}\label{prop2.3}
The FBSDE (\ref{eq4a}) admits a set of adapted solutions $(x_i,p_i, i=1,\cdots,N)$ if and only if
   (\ref{eq10}) admits an adapted solution $(x^{(N)},p^{(N)})$.
\end{proposition}
\emph{Proof.} If
   (\ref{eq10}) admits an adapted solution $(x^{(N)},p^{(N)})$, then (\ref{eq4a}) is decoupled. The existence of a set of solutions to (\ref{eq4a}) follows. The part of necessity is straightforward. \hfill$\Box$


We further discuss the optimal feedback control of (P1$^{\prime}$).
Let $\check{p}^{(N)}(t)=P(t)\check{x}^{(N)}(t)+\check{s}(t)$, $t\geq0$, where $P\in \mathbb{R}^{n\times n}$ and $\check{s}\in \mathbb{R}^n$.  Then by (\ref{eq10}) we have
$$\begin{aligned}
  d\check{p}^{(N)}=&P\Big[(A+G)\check{x}^{(N)}+B\check{u}^{(N)}-R_2^{-1}\check{p}^{(N)})dt
  \\&+\frac{1}{N}\sum_{i=1}^N\sigma dW_i\Big]+\dot{P}\check{x}^{(N)}dt+d\check{s}\cr
  =&-\Big[(A+G )^T(P\check{x}^{(N)}+\check{s})+(\Psi-Q)\check{x}^{(N)}
  \\&+\bar{\eta}\Big]dt
+\frac{1}{N}\sum_{i=1}^N\sum_{j=1}^N{\beta}_i^jdW_j.
\end{aligned}$$
This implies 
\begin{align}\label{eq11}
&\dot{P}+(A+G)^TP+P(A+G)-PR_2^{-1}P\cr
&-(I-\Gamma)^TQ(I-\Gamma)=0,\ P(T)=-H,\\
\label{eq12}
&d\check{s}+\big[(A+\bar{G})^T\check{s}+PB\check{u}^{(N)}+\bar{\eta}\big]dt
\cr&+\frac{1}{N}\sum_{i=1}^N\sum_{j=1}^N\big(\frac{\sigma}{N}-\beta_i^j
\big)dW_j
=0,\ \check{s}(T)=0,
\end{align}
where $\bar{G}\stackrel{\Delta}{=}G-R_2^{-1}P$.

 By the local Lipschitz continuous property of the quadratic (matrix) function, (\ref{eq11}) must admit a unique local solution in a small time duration $[T_0, T]$.
The global existence of the solution for $t\in [-\infty, T]$ or $[0, T]$ can be referred to \cite{FJ96, AFIJ03}. From Proposition \ref{prop1b}, we obtain that under (A2$^{\prime}$), (\ref{eq11}) has a unique solution in $C([0,T], \mathbb{R}^{n\times n})$.

\begin{theorem}\label{thm2}
  Under (A0), (A1), (A2$^{\prime}$), Problem (P1$^{\prime}$) has a minimizer 
 \begin{equation}\label{eq11b}
\check{f}(t)=-{R_2^{-1}}[P(t)\check{x}^{(N)}(t)+\check{s}(t)],\ t\geq0,
\end{equation}
where $P$ and $\check{s}$ are solutions of (\ref{eq11}) and (\ref{eq12}), respectively.
\end{theorem}
\emph{Proof. } Under (A2$^{\prime}$), (\ref{eq11}) admits a  unique solution $P$, which implies (\ref{eq12}) has a unique solution $\check{s}$ in $C([0,T],\mathbb{R}^n)$. By \cite[Theorem 2.4.1]{MY99}, (\ref{eq10}) admits a unique solution $(\check{x}^{(N)},\check{p}^{(N)} ),$ where $\check{p}^{(N)}=P\check{x}^{(N)}+\check{s}$. From Proposition \ref{prop2.3}, (\ref{eq4a}) is solvable. This with  Theorem \ref{thm1} completes the proof. \hfill{$\Box$}

\begin{remark}
      From the above analysis, (A2$^{\prime}$) is sufficient for solvability of FBSDE (\ref{eq4a}). Indeed, from \cite{MY99}, (A2$^{\prime}$) is also a necessary condition to ensure that the solvability of (\ref{eq4a}) holds for any $\check{u}_i\in {\mathcal U}_{c}^{\rm F}$.
\end{remark}

\subsection{Distributed Strategy Design}\label{sec4}

After the ``worst-case" drift $\check{f}$ is applied, we have the following optimal control
problem.

 \textbf{(P2)}: Minimize ${J}_{\rm soc}^{\rm F}(u,\check{f}(u))$ over $\{u_i\in {\mathcal U}_{c}^{\rm F}, i=1,\cdots,N\}$, where
\begin{align}\label{eq12a}
dx_i =& [Ax_i+Bu_i+Gx^{(N)}\!-\!{R_2^{-1}}(Px^{(N)}+s)]dt\cr
&+\sigma dW_i,\  1\leq i\leq N,\\ \label{eq12b}
ds=&-\big[(A+\bar{G})^Ts+PB{u}^{(N)}+\bar{\eta}\big]dt
\cr
&+\frac{1}{N}\sum_{i=1}^N\sum_{j=1}^N(\beta_i^j-\frac{\sigma}{N})dW_j
, s(T)=0.\\
\label{eq12c}
 {J}_{\rm soc}^{\rm F}(u)=&\frac{1}{2}\sum_{i=1}^N
\mathbb{E}\int_0^{T}\Big\{\big\|x_i
   -\Gamma x^{(N)}-\eta\big\|^2_{Q}+
   \|u_i\|^2_{R_1}\cr&\!-\!\|Px^{(N)}+s\|^2_{R_2^{-1}}\Big\}dt
\!+\!\frac{1}{2}\mathbb{E}\|x_i(T)\|_{H}^2.
\end{align}

We first show that Problem (P2) has the property of uniformly convexity under certain conditions.
\begin{lemma}\label{lem1a}
\textcolor{blue}{Assume that A0), A1), A2$^{\prime}$) hold. There exists a sufficiently large $C_0>0$ with $R_1>C_0I$ and $R_2>C_0I$ 
  such that Problem (P2) is uniformly convex in $\textbf{u}$.} 
\end{lemma}
\emph{Proof.} 
   Denote ${\bf{\bar{Q}}}=diag\{Q,\cdots,Q\}-\frac{1}{N}\textbf{11}^T\otimes (\Psi +PR_2^{-1}P)$, $\textbf{R}_1=diag\{R_1,\cdots,R_1\}$, ${\bf{\bar{A}}}=diag\{A,\cdots,A\}+\frac{1}{N}\textbf{11}^T\otimes (G-R_2^{-1}P)$. By 
   a similar argument with \cite{LZ99}, we obtain that Problem (P2) is uniformly convex if for any $u_i\in {\mathcal U}_{c}^{\rm F}$,
   $$
   \begin{aligned}
    & \mathbb{E}\int_0^{T}\Big(\textbf{z}^T{\bf{\bar{Q}}}\textbf{z}
  +\textbf{u}^T\textbf{R}_1\textbf{u}-N\grave{s}^{T}R_2^{-1}\grave{{s}}\Big)dt
\\&+\mathbb{E}\|\textbf{z}(T)\|_{\bf{H}}^2
\geq \delta\mathbb{E} \int_0^T\|\textbf{u}\|
^2dt  \end{aligned}$$ where  $\textbf{z}\in \mathbb{R}^{nN}$ and $\grave{s}\in \mathbb{R}^{n}$ satisfy
\begin{align}
  \label{eq13ba}
   d\textbf{z} = &(\bar{\textbf{A}}\textbf{z}+\textbf{B}\textbf{u}-\textbf{1}\otimes R_2^{-1} \grave{s} )dt,\
   \textbf{z}(0)= 0,\\ \label{eq13bb}
   d\grave{s}=&-\big[(A+\bar{G})^T\grave{s}+\frac{1}{N}PB(\textbf{1}^T\otimes I) \textbf{u}\big]dt
\cr&+\frac{1}{N}\sum_{i=1}^N\sum_{j=1}^N\grave{\beta}_i^jdW_j
,\ \grave{s}(T)=0.
\end{align}
By \cite[Chapter 7]{YZ99} and (\ref{eq13bb}),
\textcolor{blue}{\begin{align}\label{eq15}
  &\mathbb{E}\int_0^T\|\grave{s}(t)\|^2dt\leq \frac{C_1}{N^2}\mathbb{E}\int_0^T \|\textbf{1}^T\otimes I\|^2\|\textbf{u}(t)\|^2dt\cr
\leq &\frac{C_1}{N}\mathbb{E}\int_0^T \|\textbf{u}(t)\|^2dt.
\end{align}
This with (\ref{eq13ba}) leads to
$\mathbb{E}\int_0^T\|\textbf{z}\|^2dt\leq C_2\int_0^T\mathbb{E}\|\textbf{u}\|^2dt$.
Note that
$$\begin{aligned}
  \lambda_{min}({\bf{\bar{Q}}})\geq& \lambda_{min}(Q)\!-
  [\lambda_{max}(\Psi )\!+\!\lambda_{max}(PR_2^{-1}P)]\\
  \geq&
-[\lambda_{max}(\Psi )+\lambda_{max}(PR_2^{-1}P)],
\end{aligned}$$}
where $ \lambda_{min}(Q)$ and $ \lambda_{max}(Q)$ are smallest and largest eigenvalues of $Q$, respectively.
From this with (\ref{eq15}),
 there exists $\delta>0$ and $C_0>0$ such that for $R_1>C_0I$ and $R_2>C_0I$,
 $$\begin{aligned}
   &\mathbb{E}\int_0^{T}\Big(\textbf{z}^T{\bf{\bar{Q}}}\textbf{z}
 \! +\!\textbf{u}^T\textbf{R}_1\textbf{u}-N\grave{s}^{T}R_2^{-1}\grave{{s}}\Big)dt
+\mathbb{E}\|\textbf{z}(T)\|_{\bf{H}}^2\cr
&\geq
\mathbb{E}\int_0^T\Big[\textbf{u}^T(\textbf{R}_1-{C}_0I_{nN})\textbf{u}\Big]dt \geq \delta \mathbb{E}\int_0^T(\textbf{u}^T\textbf{u})dt.
 \end{aligned}$$
\hfill $\Box$

\subsubsection{The Social Variational Derivation}

Note that the social optimum implies the person-by-person optimality \cite{H80}. \textcolor{blue}{If the social cost is convex, then the socially optimal solution exists and coincides with the person-by-person optimal solution \cite{WZ17}.}
We now provide a transformation of the original social optimum problem by variational derivation and person-by-person optimality. Suppose that $\hat{u}=(\hat{u}_1,\cdots,\hat{u}_N)$ is a minimizer to
Problem (P2), where $\hat{u}_j\in {\mathcal{U}}_c^{\rm F}$. Let $\hat{x}_j$ correspond to $\hat{u}_j$, $j=1,\cdots,N$ and $\hat{x}^{(N)}=\frac{1}{N}\sum_{j=1}^N\hat{x}_j$. Let $\hat{s}$ correspond to $\hat{u}_1$, $\cdots$, $\hat{u}_N$. Fix $\hat{u}_{-i}=(\hat{u}_1,\cdots,\hat{u}_{i-1},\hat{u}_{i+1},\cdots,\hat{u}_N)$, and perturb $u_i$. Denote $\delta u_i=u_i-\hat{u}_i$, $\delta x_j=x_j-\hat{x}_j$, $\delta x^{(N)}=\frac{1}{N}\sum_{j=1}^N\delta x_j$ and $\delta s=s-\hat{s}$. Let the strategy variation $\delta {u}_i$ be adapted to ${\mathcal F}_t$
and satisfy $\mathbb{E}\int_0^{T}  \|\delta  {u}_i\|^2dt<\infty.$ Let $\delta {J}_i$ be the variation of ${J}_i$ with respect to $\delta u_i$.
By (\ref{eq12a}) and (\ref{eq12b}),
\begin{align}
  \frac{d\delta x_j}{dt}=&A\delta x_j+\frac{\bar{G}}{N}\delta x_i+\frac{\bar{G}}{N}\sum_{k\not=i}\delta x_k-R_2^{-1}\delta s,\cr
  & \ j\not =i,\ \delta x_j(0)=0,\cr
  d\delta s=&-\big[(A+\bar{G})^T\delta s+\frac{1}{N}PB\delta u_i\big]dt\cr
&+\frac{1}{N}\sum_{i=1}^N\sum_{j=1}^N\delta \beta_i^jdW_j
,\ \delta s(T)=0.\label{eq20aa}
\end{align}
where $\bar{G}\stackrel{\Delta}{=}G-R_2^{-1}P$. This implies $\delta x_j=\delta x_{k}$,\ for any $j,k\not =i$. Thus,
$$ \frac{d\delta x_j}{dt}\!=\!(A+\frac{N-1}{N}\bar{G})\delta x_j+\frac{\bar{G}}{N}\delta x_i\!-\!R_2^{-1}\delta s, \delta x_j(0)\!=\!0,$$
which gives that
$$\delta x_j(t)\!=\!\!\int_0^t\!\!e^{(A+\frac{N-1}{N}\bar{G})(t-\tau)}\big(\frac{\bar{G}}{N}\delta x_i(\tau)-R_2^{-1}\delta s(\tau)\big)d\tau\!.$$
We further have
$$\begin{aligned}
  \delta x^{(N)}(t)=&\frac{1}{N}\delta x_i(t)+\frac{N-1}{N}\int_0^te^{(A+\frac{N-1}{N}\bar{G})(t-\tau)}\\
&\times\Big(\frac{\bar{G}}{N}\delta x_i(\tau)-R_2^{-1}\delta s(\tau)\Big)d\tau.
\end{aligned}$$
By this with (\ref{eq12b}), one can obtain
\begin{align*}
  \delta J_i(u,\check{f})=&\mathbb{E}\int_0^T  \Big\{[\hat{x}_i-\Gamma \hat{x}^{(N)}-\eta]^T Q\big[\delta x_i-{\Gamma}\delta x^{(N)}\big]\cr
  &-(P\hat{x}^{(N)}+s)^TR_2^{-1}(P\delta x^{(N)}+\delta s)\cr
  &+\hat{u}_i^TR_1\delta u_i\Big\}dt+\mathbb{E}[x_i^T(T)H\delta x_i(T)],
  \end{align*}
and for $j\not = i$,
  \begin{align*}
  &\delta J_j(u,\check{f})\cr
  =&\mathbb{E}\int_0^T   \big[(\hat{x}_j-\Gamma \hat{x}^{(N)}-\eta)^T Q(\delta x_j-\Gamma \delta x^{(N)})\cr
  &-(P\hat{x}^{(N)}+s)^TR_2^{-1}(P\delta x^{(N)}+\delta s)\big]dt\cr
  &+ \mathbb{E}[x_j^T(T)H\delta x_j(T)].
    \end{align*}
The above equation further implies that
 \begin{align*}
&\sum_{j\not= i}\delta J_j(u,\check{f})\cr
  =&\mathbb{E}\int_0^T   \big(\hat{x}_{-i}^{(N)}-\textstyle{\frac{N-1}{N}}(\Gamma \hat{x}^{(N)}+\eta)\big)^T
 Q\Big[(I-\textstyle{\frac{N-1}{N}}\Gamma) \cr&\cdot\int_0^te^{(A+\frac{N-1}{N}\bar{G})(t-\tau)} (\bar{G}\delta x_i-NR_2^{-1}\delta s)d\tau-{\Gamma}\delta x_i\Big]\cr
 &-(P\hat{x}^{(N)}+s)^TR_2^{-1}\Big[P\big(\textstyle{\frac{(N-1)^2}{N^2}}
 \int_0^te^{(A\!+\!\frac{N-1}{N}\bar{G})(t-s)} \cr
 &\cdot\bar{G}\delta x_ids\!+\!\textstyle{\frac{N-1}{N}}\delta x_i\big)\!+\!(N-1)\delta s\Big]dt
+\mathbb{E}\Big[x_{-i}^{(N)}(T)^T \cr&\cdot H\int_0^Te^{(A+\frac{N-1}{N}\bar{G})(T-t)}(\bar{G}\delta x_i-NR_2^{-1}\delta s)dt\Big],
     \end{align*}
    where $\hat{x}_{-i}^{(N)}=\frac{1}{N}\sum_{j\not =i}\hat{x}_j$.
    Let $\delta\psi_i=N\delta s$. Then from (\ref{eq20aa}),
  \begin{align}\label{eq20a}
    d\delta\psi_i=&-\big[(A+\bar{G})^T\delta \psi_i+PB\delta u_i\big]dt\cr
&+\sum_{i=1}^N\sum_{j=1}^N\delta \beta_i^jdW_j
, \delta \psi_i(T)=0.
  \end{align}
      \textcolor{blue}{For large $N$, it is plausible to approximate $\hat{x}^{(N)}$ by a deterministic function $\bar{x}$.} The zero first-order variational condition combined with the mean field approximation gives
 \begin{align}\label{eq18}
 &\mathbb{E}\int_0^T  \Big\{(\hat{x}_i-\Gamma \bar{x}-\eta)^T
 Q\delta x_i-  \big[((I-\Gamma) \bar{x}-\eta)^TQ\Gamma\cr
 &+(P\bar{x}+\bar{s})^TR_2^{-1}P\big]\delta x_i  +\big[ ((I-\Gamma) \bar{x}-\eta)^TQ(I-\Gamma)\cr
  &\!-\!(P\bar{x}\!+\!\bar{s})^TR_2^{-1}P\big]
 \!\! \int_0^te^{(A+\bar{G})(t-\tau)} (\bar{G}\delta x_i\!-\!R_2^{-1}\delta \psi_i)d\tau
  \cr&-(P\bar{x}+\bar{s})^TR_2^{-1}\delta \psi_i+\hat{u}_i^TR_1\delta u_i\cr&+\bar{x}^T(T)He^{(A+\bar{G})(T-t)}(\bar{G}\delta x_i-R_2^{-1}\delta \psi_i)\Big\}dt\cr
  &+\mathbb{E}[x_i^T(T)H\delta x_i(T)]=0,
  \end{align}
where $\bar{x}\in C([0,T], \mathbb{R}^n)$ is an approximation of $\hat{x}^{(N)}$. From observation, the equation (\ref{eq18}) is the zero variation condition for the optimal control problem
with the cost function:
 \begin{align}\label{eq17}
&J_i^{\prime}(u_i)\cr
=&\frac{1}{2}\mathbb{E}\int_0^T \Big\{x^T_iQx_i+2\big[- \left((I-\Gamma) \bar{x}-\eta\right)^TQ\Gamma\cr
&-(\Gamma \bar{x}+\eta)^T
   Q -(P\bar{x}+\bar{s})^TR_2^{-1}P\big]x_i
 \cr&    +2\big[ ((I-\Gamma) \bar{x}-\eta)^TQ(I-\Gamma)-(P\bar{x}+\bar{s})^TR_2^{-1}P\big]
 \cr&\times \int_0^te^{(A+\bar{G})(t-\tau)} (\bar{G} x_i-R_2^{-1} \psi_i)d\tau
  \cr&+\bar{x}^T(T)He^{(A+\bar{G})(T-t)}(\bar{G} x_i-R_2^{-1} \psi_i)\cr
  &+{u}_i^TR_1 u_i\Big\}dt+\frac{1}{2}\mathbb{E}[x_i^T(T)H x_i(T)]\cr
=& \frac{1}{2}\mathbb{E}\int_0^T\!\!\Big\{x_i^TQ x_i-2\big[\Psi \bar{x}+\bar{\eta} +(P\bar{x}+\bar{s})^TR_2^{-1}P\big]x_i
 \cr&   +2 v^T(\bar{G} x_i-R_2^{-1} \psi_i) 
 +u_i^TR_1u_i \Big\}dt\cr
 &+\frac{1}{2}\mathbb{E}[x_i^T(T)H x_i(T)],
  \end{align}
  where the second equality holds by an exchange of order of
the integration, and
\begin{align*}
v(t)\stackrel{\Delta}{=}&\int_t^{T}e^{({A+\bar{G}})^T(\tau-t)}\big[ (I-\Gamma)^TQ((I-\Gamma) \bar{x}-\eta)\cr
&-PR_2^{-1}(P\bar{x}+\bar{s})\big]d\tau+ e^{(A+\bar{G})^T(T-t)}H\bar{x}(T).
\end{align*}

\subsubsection{Mean field approximation}\label{sec4a}

 Based on (\ref{eq12a}), (\ref{eq20a}) and (\ref{eq17}), we construct the following auxiliary optimal control problem.

 \textbf{(P3)}: minimize $ \bar{J}_i(u_i)$ over $u_i\in{\cal U}_i$, where
\begin{align}\label{eq19}
d\grave{x}_i =& [A\grave{x}_i+Bu_i+G\bar{x}-{R_2^{-1}}(P\bar{x}+\bar{s})]dt\cr
&+\sigma dW_i,\  \grave{x}_i(0)=x_{i0},\\
d\psi_i=&-\big[(A+\bar{G})^T \psi_i+PB u_i\big]dt
\cr
&+ z_i^idW_i
,\  \psi_i(T)=0.\label{eq19c}\\
\label{eq20}
 \bar{J}_i(u_i)=& \frac{1}{2}\mathbb{E}\int_0^T  \Big[\grave{x}_i^TQ \grave{x}_i-2(P\bar{x}+\bar{s})^TR_2^{-1}P\grave{x}_i
 \cr
 &-2(\Psi \bar{x}+\bar{\eta})^T \grave{x}_i +2v^T(\bar{G} \grave{x}_i-R_2^{-1} \psi_i)
 \cr& -2(P\bar{x}+\bar{s})^TR_2^{-1} \psi_i+u_i^TR_1u_i \Big]dt\cr
 &+\frac{1}{2}\mathbb{E}[\grave{x}_i^T(T)H \grave{x}_i(T)].
\end{align}
Here $\bar{s},v$ are determined by
\begin{align} \label{eq25a}
\dot{\bar{s}}&=-(A+G-R_2^{-1}P)^T\bar{s}-PB\bar{u}-\bar{\eta},\ \bar{s}(T)=0,\\
 \label{eq25}
    \dot{v}&=-(A+\bar{G})^Tv-\big[(I-\Gamma)^TQ((I-\Gamma)\bar{x}-\eta)\cr
    &-PR_2^{-1}(P\bar{x}+\bar{s})\big], \  v(T)=H\bar{x}(T),
  \end{align}
and  $\bar{x},\bar{u}\in C([0,T], \mathbb{R}^n)$ are given for approximations to $\hat{x}^{(N)},\hat{u}^{(N)}$, respectively.

  \begin{theorem}\label{thm3}
Assume that (A0)-(A1), (A2$^{\prime}$) hold. Problem (P3) has a unique optimal control
\begin{equation}\label{eq24}
  \hat{u}_i(t)=-R_1^{-1}B^T[k_i(t)-P(t)l(t)],\   1\leq i\leq N, 
\end{equation}
where $(l,k_i, \zeta_i)$ is a unique adaptive solution to the following (decoupled) FBSDE
\begin{align}\label{eq23}
dl=&[(A+\bar{G})l+R_2^{-1}v+R_2^{-1}(P\bar{x}+\bar{s})]dt, l(0)=0,\\
 dk_i=&-\Big\{A^T k_i+Q\grave{x}_i-(\Psi \bar{x}+\bar{\eta})-PR_2^{-1}(P\bar{x}+\bar{s})\cr
 &+\bar{G}^Tv\Big\}dt+\zeta_idW_i,\  k_i(T)=H\grave{x}_i(T).\label{eq23b}
\end{align}
\end{theorem}
\emph{Proof.} Since $Q\geq0$ and $R_1>0$, then from \cite{LZ99,HH16}, (P3) is uniformly convex in $u_i$ and there exists a unique optimal control for (P3),
denoted as $\hat{u}_i$. Then
 \begin{align}\label{eq25a}
0=&\delta {\bar{J}}_i(\hat{u}_i)\cr
=&\mathbb{E}\int_0^T \! \Big[(Q\grave{x}_i \! - \! \Psi \bar{x} \! - \! \bar{\eta})^T\delta \grave{x}_i
-(P\bar{x}+\bar{s})^TR_2^{-1}P\delta \grave{x}_i \cr& +v^T(\bar{G} \delta \grave{x}_i-R_2^{-1}\psi_i)-(P\bar{x}+\bar{s})R_2^{-1}\psi_i\cr
 &+u_i^TR_1\delta u_i \Big]dt +\mathbb{E}[\grave{x}_i^T(T)H\delta \grave{x}_i(T)],
\end{align}
where $\delta u_i=u_i-\hat{u}_i$, $\delta \grave{x}_i=\grave{x}_i-\hat{x}_i$, and $\delta \psi_i=\psi_i-\hat{\psi}_i$.
Note that (\ref{eq23}) and (\ref{eq23b}) are decoupled. Given $\bar{x}, \bar{u}\in C([0,T], \mathbb{R}^n)$, (\ref{eq23}) is a standard linear BSDE and so has a unique solution $(k_i,\zeta_i)$.
 Note that
\begin{align*}
d(\delta \grave{x}_i)=&(A\delta \grave{x}_i+B\delta u_i)dt,\cr
    d(\delta\psi_i)=&-\big[(A+\bar{G})^T\delta \psi_i+PB\delta u_i\big]dt
\cr
&+\sum_{i=1}^N\sum_{j=1}^N\delta \beta_i^jdW_j
,\ \delta \psi_i(T)=0.
\end{align*}
By It\^{o}'s formula, we have
$$  \begin{aligned}
&\mathbb{E}[\grave{x}_i^T(T)H\delta \grave{x}_i(T)]\cr
=&  \mathbb{E}[k_i^T(T)\delta \grave{x}_i(T)- k_i^T(0)\delta \grave{x}_i(0)]\cr
=&\mathbb{E}\int_{0}^{T} \Big\{-\big[ Q\grave{x}_i-(\Psi \bar{x}+\bar{\eta}) \\
 &-PR_2^{-1}(P\bar{x}+\bar{s})+\bar{G}^Tv\big]^T\delta \grave{x}_i+ k_i^TB\delta u_i\Big\}dt,
\end{aligned}$$
and $$
\begin{aligned}
  &\mathbb{E}[l^T(T)\delta \psi_i(T)- l^T(0)\delta \psi_i(0)]\cr
=&\mathbb{E}\!\int_0^T\!\!\big[(R_2^{-1}v+R_2^{-1}(P\bar{x}+\bar{s}))^T\delta \psi_i\!-\!l^TPB\delta u_i\big]dt.
\end{aligned}$$
This and (\ref{eq25a}) gives
$$0=\mathbb{E}\int_{0}^{T} (R_1u_i+B^Tk_i-B^TPl)^T\delta u_idt,$$
which implies $ \hat{u}_i=R_1^{-1}B^T(Pl-k_i), \ 1\leq i\leq N$.  \hfill $\Box$

Let $k_i=K\grave{x}_i+\varphi$. Then by (\ref{eq19}) and (\ref{eq23}),
\begin{align*}
  dk_i=&K(A\grave{x}_i-BR_1^{-1}B^T(K\grave{x}_i-Pl+\varphi)+G\bar{x}\cr
  &-{R_2^{-1}}(P\bar{x}+\bar{s}))dt+K\sigma dW_i+\dot{K}\grave{x}_i+\dot{\varphi}\cr
  =&-\Big\{A^T (K\grave{x}_i+\varphi)+Q\grave{x}_i-(\Psi \bar{x}+\bar{\eta})
 \cr
&-PR_2^{-1}(P\bar{x}+\bar{s})+\bar{G}^Tv\Big\}+\zeta_idW_i,
  \end{align*}
which implies
\begin{align}\label{eq26}
&\dot{K}+A^TK+KA-KBR_1^{-1}B^TK+Q=0,\ \cr& \  K(T)=H,\\
  \label{eq27}
&\dot{\varphi}+(A-BR_1^{-1}B^TK)^T\varphi+KBR_1B^TPl
\cr&+K\bar{G}\bar{x}-KR^{-1}_2\bar{s}
  -(\Psi \bar{x}+\bar{\eta})\cr&
 -PR_2^{-1}(P\bar{x}+\bar{s})+\bar{G}^Tv=0, \ \varphi(T)=0.
\end{align}
Besides, applying (\ref{eq24}) into (\ref{eq19}), we obtain
\begin{align*}
d{\hat{x}}^{(N)} =& [A{\hat{x}}^{(N)}-BR_1^{-1}B^T(K{\hat{x}}^{(N)}-Pl+\varphi)\cr
&+\bar{G}{\bar{x}}-{R_2^{-1}}\bar{s}]dt+\frac{1}{N}\sum_{i=1}^N\sigma dW_i,\quad
\end{align*}
where ${\hat{x}}^{(N)}(0)=\frac{1}{N}\sum_{i=1}^N{x}_{i0}.$ As an approximation, one can obtain
\begin{align}\label{eq28}
  \dot{\bar{x}} =& (\bar{A}\!+\!\bar{G})\bar{x}+BR_1^{-1}B^T(Pl\!-\!\varphi)\!-\!{R_2^{-1}}\bar{s}, \bar{x}(0)=\bar{x}_0,
\end{align}
where $\bar{A}\stackrel{\Delta}{=}A-BR_1^{-1}B^TK$.

By (\ref{eq28}), (\ref{eq20}), (\ref{eq25}), and (\ref{eq27}),
\begin{equation}\label{eq29}
  \left\{
\begin{aligned}
\dot{\bar{x}} =& (\bar{A}+\bar{G})\bar{x}+BR_1^{-1}B^T(Pl-\varphi)\cr&-{R_2^{-1}}\bar{s},
 \  \bar{x}(0)=\bar{x}_0,\cr
\dot{l}=&(A+\bar{G})l+R_2^{-1}v+R_2^{-1}(P\bar{x}+\bar{s}),\ l(0)=0,\cr
\dot{\bar{s}}=&-(A+\bar{G})^T\bar{s}+PBR_1^{-1}B^T(K\bar{x}-P l+\varphi)\cr&-\bar{\eta},\ \bar{s}(T)=0, \cr
\dot{\varphi}=&-\bar{A}^T\varphi-KBR_1^{-1}B^TPl-K\bar{G}\bar{x}+KR^{-1}_2\bar{s}\cr
&+\Psi \bar{x}+\bar{\eta}
 +PR_2^{-1}(P\bar{x}+\bar{s})\!-\!\bar{G}^Tv,  \varphi(T)=0,\cr
 \dot{v}=&-(A+\bar{G})^Tv+(\Psi-Q)\bar{x}+\bar{\eta}\cr
 &+PR_2^{-1}(P\bar{x}+\bar{s}), \ v(T)=H\bar{x}(T).
\end{aligned}
\right.
\end{equation}

 For further analysis, we assume:

(\textbf{A3}) (\ref{eq29}) admits a unique solution in $C([0,T], \mathbb{R}^{5n})$.

Note that (\ref{eq29}) can be taken as an FBSDE without diffusion terms. The condition of contraction mapping in Theorem 5.1 of \cite{MY99} holds necessarily.
Thus, (\ref{eq29}) must admit a unique solution in a small time duration $[T_0,T]$.
However, some additional conditions are needed for existence of a (global) solution to (\ref{eq29}) in the time duration $[0,T]$.
We now give a sufficient condition that ensures ({A3}).

Let
$$M_{11}=\left[\begin{array}{cc}
  \bar{A}+\bar{G}& BR_1^{-1}B^TP\\
  R_2^{-1}P& A+\bar{G}
\end{array}
\right],$$
$$M_{12}=\left[\begin{array}{ccc}
 -R_2^{-1}& -BR_1^{-1}B^T& 0\\
 R_2^{-1}&0& R_2^{-1}
 \end{array}
\right],$$
$$M_{21}=\left[\begin{array}{cc}
  PBR_1^{-1}B^TK&\!-PBR_1^{-1}B^TP\!\\
  -K\bar{G}+\Psi+PR_2^{-1}P&0\\
  \Psi-Q+PR_2^{-1}P&0
\end{array}\right]\!\!,$$
$$M_{22}=
\left[\begin{array}{ccc}
 -(A+\bar{G})&  \!\!PBR_1^{-1}B^TK\!\!& 0\\
  \!(K+P)R_2^{-1}\!&-\bar{A}&-\bar{G}^T\\
 PR_2^{-1}&0&\!-(A+\bar{G})\!
\end{array}\right]\!\!.
$$
Then (\ref{eq29}) can be written as
\begin{equation}\label{eq26b}
\left[\begin{array}{c}
 \dot{\bar{x}}\\
 \dot{{l}}\\
 \dot{s}\\
\dot{ \varphi}\\
\dot{v}
\end{array}\right]=\left[\begin{array}{cc}
M_{11}& M_{12}\\
M_{21}&M_{22}
\end{array}\right]\left[\begin{array}{c}
 {\bar{x}}\\
  l\\
 {s}\\
{ \varphi}\\
{v}
\end{array}\right]+\left[\begin{array}{c}
0\\
0\\
-\bar{\eta}\\
\bar{\eta}\\
\bar{\eta}
\end{array}\right].
\end{equation}
\begin{proposition}\label{prop4a}
If the Riccati differential equation
$$\begin{aligned}
  \dot{Y}&=M_{21}+M_{22}Y-YM_{11}-YM_{12}Y,\\
 &\textstyle{ Y(T)=\left[
\begin{array}{ccc}
0&\ 0&\ H^T\\
0&\ 0&\ 0\end{array}
\right]^T}
\end{aligned}$$
admits a solution $Y\in \mathbb{R}^{3n\times 2n}$ in $[0,T]$, then ({A3}) holds.
  Furthermore, under the assumption $\bar{\eta}=0$,
if the Riccati differential equation
\begin{equation}\label{eq38a}
  \dot{Z}=M_{12}+M_{11}Z-ZM_{22}-ZM_{21}Z,  Z(0)=0
\end{equation}
admits a solution $Z\in \mathbb{R}^{2n\times 3n}$ in $[0,T]$, then ({A3}) holds.
\end{proposition}
\emph{Proof.} Denote $m = [\bar{x}^T,l^T]^T$, and $z=[s^T,\varphi^T,v^T]^T$. Let $z=Y m+\alpha$, $\alpha(T)=0$. Then $\textstyle{ Y(T)=\left[
\begin{array}{ccc}
0&\ 0&\ H^T\\
0&\ 0&\ 0\end{array}
\right]^T}\!.$ By (\ref{eq26b}),
\begin{align*}
  \dot{z}=&\dot{Y}m+Y(M_{11}m+M_{12}z)+\dot{\alpha}\cr
  =& (\dot{Y}+YM_{11}+YM_{12}Y)\bar{x}+YM_{12}\alpha+\dot{\alpha}\cr
  =&M_{21}\bar{x}+M_{22}(Y\bar{x}+\alpha)+[-\bar{\eta}^T,\bar{\eta}^T,\bar{\eta}^T]^T.
\end{align*}
Thus, we obtain
\begin{align}\label{eq27a}
\dot{Y}=&M_{21}+M_{22}Y-YM_{11}-YM_{12}Y, \\
\label{eq27b}
\dot{\alpha}=&(M_{22}-YM_{12})\alpha+[-\bar{\eta}^T,\bar{\eta}^T,\bar{\eta}^T]^T,
\end{align}
where $\textstyle{ Y(T)=\left[
\begin{array}{ccc}
0&\ 0&\ H^T\\
0&\ 0&\ 0\end{array}
\right]^T}$ and $\alpha(T)=0$. Since (\ref{eq27a}) admits a solution, then (\ref{eq27b}) has a solution. Applying $z=Ym+\alpha$ into (\ref{eq26b}), we have
$$
\dot{m}=M_{11}m+M_{12}(Ym+\alpha),\ m(0)=[\bar{x}_0^T,0]^T,
$$
which implies (\ref{eq29}) admits a unique solution in $[0,T]$.

Denote $Z=\textstyle{\left[\begin{array}{ccc}
Z_{11}&Z_{12}&Z_{13}\\
Z_{21}&Z_{22}&Z_{23}
\end{array}\right]}$. Note that ${s}(T)=\varphi(T)=0$, and $v(T)=H\bar{x}(T)$.
We have
$$\begin{aligned}
&\left[\begin{array}{cc}
0&0\\0&
0\\H&0
\end{array}\right]
\left[\begin{array}{ccc}
Z_{11}&Z_{12}&Z_{13}\\
Z_{21}&Z_{22}&Z_{23}
\end{array}\right]
\\
=&\left[\begin{array}{ccc}
0&0&0\\
0&0&0\\
HZ_{11}(T)&HZ_{12}(T)&HZ_{13}(T)
\end{array}\right]\not =I_{3n}.
\end{aligned}$$
By the modified Radon's Lemma (see e.g., \cite[Theorem 3.1.3]{{AFIJ03}}), the proposition follows.
  \hfill $\Box$

\subsection{Asymptotic Optimality}

For Problem (PF), we may design the following decentralized control:
\begin{equation}
\label{eq24b}
  \hat{u}_i(t)=-R^{-1}_{1}B^T[K(t)x_i(t)-P(t)l(t)+\varphi(t)],
  \end{equation}
where $K,P$ are given by (\ref{eq26}) and (\ref{eq11}), respectively, and $l$ and $\varphi$ are determined by (\ref{eq29}).
  After the control laws (\ref{eq24b}) are applied, we obtain the state equations of agents as follows:
 \begin{align}\label{eq24c}
 d\hat{x}_i=&[\bar{A}\hat{x}_i+\bar{G}\hat x^{(N)}\!+\!BR_1^{-1}B^T(Pl-\varphi)-{R_2^{-1}}\hat{s}]dt\cr
 &+\sigma dW_i,\ i=1,\cdots,N.
 \\ \label{eq24d}
 d\hat{s}=&-\big[(A+\bar{G})^T\hat{s}+PB\hat{u}^{(N)}+\bar{\eta}\big]dt\cr
&+\frac{1}{N}\sum_{i=1}^N\sum_{j=1}^N(\beta_i^j-\frac{\sigma}{N})dW_j
,\ \hat{s}(T)=0.
\end{align}

For further analysis, we assume

(\textbf{A4)} The Riccati equation admits a solution:
\begin{equation*}
\begin{split}
  &\dot{\tilde{P}}+\tilde{P}(\bar{A}+\bar{G})+({A}+\bar{G})^T\tilde{P}-\tilde{P}R_2^{-1}\tilde{P}\cr
  &+PBR_1^{-1}B^TK=0,\ \tilde{P}(T)=0.
\end{split}
\end{equation*}

\begin{lemma}\label{lem1}
Assume that (A0)-(A1), (A2$^{\prime}$), (A3)-(A4) hold. For the system (\ref{eq1})-(\ref{eq2}), we have 
\begin{equation}\label{eq30a}
  \sup_{0\leq t\leq T}\mathbb{E}\big(\|\hat{x}^{(N)}-\bar{x}\|^2+\|\hat{s}-\bar{s}\|^2\big)=O(1/N).
  \end{equation}
\end{lemma}
\emph{Proof.}  It follows by (\ref{eq24c}) and (\ref{eq24d}) that
\begin{equation*}
\begin{aligned}
 d\hat{x}^{(N)}=&[(\bar{A}+\bar{G})\hat x^{(N)}+BR_1^{-1}B^T(Pl-\varphi)\cr
 &-{R_2^{-1}}\hat{s}]dt+\frac{1}{N}\sum_{i=1}^N\sigma dW_i, \\
d\hat{s}=&(A+\bar{G})^T\hat{s}+PB\hat{u}^{(N)}+\bar{\eta}\big]dt\cr
&+\frac{1}{N}\sum_{i=1}^N\sum_{j=1}^N(\beta_i^j-\frac{\sigma}{N})dW_j
,\ \hat{s}(T)=0.
\end{aligned}
\end{equation*}
Denote $\xi\stackrel{\Delta}{=}\hat{x}^{(N)}-\bar x$ and $\chi\stackrel{\Delta}{=}\hat{s}-\bar{s}$. From the above two equations and (\ref{eq29}),
\begin{align}\label{eq45a}
  d\xi=&(\bar{A}+\bar{G})\xi dt-{R_2^{-1}}\chi dt+\frac{1}{N}\sum_{i=1}^N\sigma dW_i,\cr
   &\quad\xi(0)=\textstyle{\frac{1}{N}\sum_{i=1}^Nx_{i0}-\bar{x}_0},\\
  \label{eq45b} d\chi=&-\big[(A+\bar{G})^T\chi+PBR_1^{-1}B^TK\xi\big]dt\cr
&+\frac{1}{N}\sum_{i=1}^N\sum_{j=1}^N(\beta_i^j-\frac{\sigma}{N})dW_j
,\ \chi(T)=0.
\end{align}
Let $\chi(t)=\tilde{P}(t)\xi(t)+\psi(t)$, $t\geq 0$. By It\^{o} formula,
$$\begin{aligned}
  d\chi=&\dot{\tilde{P}}\xi+\tilde{P}\Big\{[(\bar{A}+\bar{G})\xi-R_2^{-1}(\tilde{P}\xi+\psi)] dt\cr
  &+\frac{1}{N}\sum_{j=1}^N\sigma dW_j
\Big\}+\dot{\psi}\\
=&-\big[(A+\bar{G})^T(\tilde{P}\xi+\psi)+PBR_1^{-1}B^TK\xi\big]dt\cr
&+\frac{1}{N}\sum_{i=1}^N\sum_{j=1}^N(\beta_i^j-\frac{\sigma}{N})dW_j,
\end{aligned}$$
which gives $\sum_{i=1}^N\beta_i^j=(\tilde{P}+I)\sigma$, and
\begin{equation*}
\begin{aligned}
&    \dot{\tilde{P}}+\tilde{P}(\bar{A}+\bar{G})+({A}+\bar{G})^T\tilde{P}-\tilde{P}R_2^{-1}\tilde{P}\cr
&+PBR_1^{-1}B^TK=0,\ \tilde{P}(T)=0,\\
 & \dot{\psi}+(A+\bar{G}-R_2^{-1}\tilde{P})^T\psi=0, \ \psi(T)=0.
  \end{aligned}
\end{equation*}
From (A4), we have $\tilde{P}$ is existent and $\psi(t)\equiv0$. Thus,
$$\begin{aligned}
\xi(t)=&e^{\Upsilon t}\xi(0)+\frac{1}{N}\int_0^te^{\Upsilon(t-\mu)}\sum_{i=1}^N
\sigma dW_i(\mu),
\end{aligned}$$
where $\Upsilon=\bar{A}+\bar{G}-R_2^{-1}\tilde{P}$. By (A0), one can obtain
$$\begin{aligned}
  \mathbb{E}\|\xi(t)\|^2
  \leq&\frac{2}{N}\big\|e^{\Upsilon t}\big\|^2\Big\{\max_{1\leq i \leq N}\mathbb{E}\|x_{i0}\|^2\!+\!\int_0^t\!\!\big\|e^{-\Upsilon\mu}\sigma\big\|^2 d\mu\Big\},
\end{aligned}$$
which completes the proof.    \hfill $\Box$

\begin{lemma}\label{lem2}
  If $u=(u_1,\cdots,u_N)$ satisfies $$\sup_{f\in \mathcal{U}_c^{\rm F}}J_{\rm soc}^{(N)}(u)\leq C,$$ then there exists $C_1$ independent of $N$ such that
  $  \mathbb{E}\int_0^T\|{u}_i\|^2dt\leq C_1$ for all $i=1,\cdots,N$.
\end{lemma}

\emph{Proof.} Let $f=0$. Since $R_1>0$, then $\sup_{f}J_{\rm soc}^{(N)}(u)\leq C$ implies
$  \mathbb{E}\int_0^T\|{u}_i\|^2dt\leq C_1$ for all $i=1,\cdots,N$. \hfill $\Box$

\begin{lemma}\label{lem4}
Assume that (A0)-(A1), (A2$^{\prime}$), (A3)-(A4) hold. Then there exists a constant $C_0$ independent of $N$ such that
$$\sup_{f\in \mathcal{U}_c^{\rm F}}\sum_{i=1}^{N} J_i(\hat{u},f)\leq NC_0.$$
\end{lemma}

\emph{Proof.} Under (A0), (A1) and (A2$^{\prime}$), $\check{f}=-R_2^{-1}(P\check{x}^{(N)}+s)$ is a maximizer of $\sum_{i=1}^{N} J_i(\hat{u},f)$, i.e., $\sum_{i=1}^{N}J_i(\hat{u},\check{f})=\sup_{f}\sum_{i=1}^{N}J_i(\hat{u},f)$. 
By Lemma \ref{lem1}, we obtain
$$\begin{aligned}
 &\sup_{0\leq t\leq T}\mathbb{E}\|\hat{x}^{(N)}\|^2\sup_{0\leq t\leq T}\big(2\mathbb{E}\|\hat{x}^{(N)}-\bar{x}\|^2+2\|\bar{x}\|^2\big)\leq C,
\end{aligned}$$
$$\sup_{0\leq t\leq T}\mathbb{E}\|s(t)\|^2\leq \sup_{0\leq t\leq T}\big(2\mathbb{E}\|\bar{s}\|^2+2\mathbb{E}\|s-\bar{s}\|^2\big)\leq C.$$
Denote $g\stackrel{\Delta}{=}\bar{G}\hat x^{(N)}+BR_1^{-1}B^T(P l-\varphi)-{R_2^{-1}}s$. Note that $l,\varphi \in C([0, T], \mathbb{R}^n)$. Then we have
$ \sup_{0\leq t\leq T} \mathbb{E}\|g(t)\|^2\leq C.$ It follows from (\ref{eq24c}) that
\begin{align*}
  \mathbb{E}\|\hat{x}_i\|^2
  \leq& C+3C\Big\|\mathbb{E}\int_{0}^Te^{\bar{A}(T-\tau)} d\tau\Big\|^2\cr
  &+3\mathbb{E}\int_{0}^ T\Big\|e^{\bar{A}(T-\tau)}\sigma\Big\|^2 d\tau.
\end{align*}
  From this with (\ref{eq24b}),
we have $$\begin{aligned}
\sum_{i=1}^{N} \sup_{f}J_i(\hat{u},f)=&\frac{1}{2}\sum_{i=1}^N
\mathbb{E}\!\int_0^{T}\!\Big\{\big\|\hat{x}_i-\Gamma \hat{x}^{(N)}-\eta\big\|^2_{Q}\cr&+\|\hat{u}_i\|^2_{R_1}-\|\check{f}\|^2_{R_2}\Big\}d t\leq NC_0.\end{aligned}$$
\hfill $\Box$

Let $\hat{k}_i\stackrel{\Delta}{=}K\hat{x}_i+\varphi$, where $\varphi$ is given by (\ref{eq29}).
We have the following approximation result. 

\begin{lemma}\label{lem5}
Assume that (A0)-(A1), (A2$^{\prime}$), (A3)-(A4) hold. Then for problem (PF), we have 
\begin{equation*}
  \sup_{0\leq t\leq T}\mathbb{E}\|\hat{k}^{(N)}-v\|^2=O(1/N),
  \end{equation*}
  where $\hat{k}^{(N)}=\frac{1}{N}\sum_{i=1}^N\hat{k}_i$ and $v$ is given by (\ref{eq29}).
\end{lemma}
\emph{Proof.} Let $\vartheta=v-K\bar{x}-\varphi$.
By (\ref{eq29}) and some elementary calculations, we obtain
$$d\vartheta(t)=-A^T\vartheta(t) dt, \quad \vartheta(T)=0,$$
which implies $\vartheta(t)\equiv 0$. This further gives $v=K\bar{x}+\varphi$.
By Lemma \ref{lem1}, we have
\begin{equation*}
   \begin{split}
  &\sup_{0\leq t\leq T}\mathbb{E}\|\hat{k}^{(N)}-v\|^2=\sup_{0\leq t\leq T}\mathbb{E}\|K(\hat{x}^{(N)}-\bar{x})\|^2\cr
  &\leq C \sup_{0\leq t\leq T}\mathbb{E}\|\hat{x}^{(N)}-\bar{x}\|^2=O(1/N).
  \end{split}
  \end{equation*}
This completes the proof.  \hfill$\Box$

We are in a position to state the result of asymptotic optimality of the decentralized control.
\begin{theorem}\label{thm4}
 Let (A0)-(A1), (A2$^{\prime}$), (A3)-(A4) hold. Assume that (P2) is convex. For Problem (PF), the set of control laws $\hat{u}=(\hat{u}_1,\cdots,\hat{u}_N)$ given by (\ref{eq24b}) has asymptotic robust social optimality, i.e.,
$$\left|\frac{1}{N} J_{\rm soc}^{\rm wo}(\hat{u})-\frac{1}{N}\inf_{u\in{\mathcal U}_c ^{\rm F}}
J_{\rm soc}^{\rm wo}({u}) \right|=O(\frac{1}{\sqrt{N}}).$$
\end{theorem}
\emph{Proof.} See Appendix A.
\hfill$\Box$

\section{Robust mean field social control over an infinite horizon}

In this section, we consider social optimum control in robust mean field model over an infinite horizon.
Let $$
\begin{aligned}
{\cal U}_{i} &=\Big\{u_i\ |\
  u_i(t)\in \sigma(x_i(s),0\leq s\leq t),\cr
  & \mathbb{E}\int_0^\infty e^{-\rho t}\|x_i(t)\|^2dt<\infty \Big\},
\end{aligned}
$$ and
\begin{eqnarray}\label{eq2b}
J_i(u,f)&\hspace*{-0.2cm}=\hspace*{-0.2cm}&
\frac{1}{2}\mathbb{E}\int_0^{\infty}e^{-\rho t}\Big\{\big\|x_i(t)
   -\Gamma x^{(N)}(t)-\eta\big\|^2_{Q}\cr
&\hspace*{-0.2cm}&+\|u_i(t)\|^2_{R_1}-\|f(t)\|^2_{R_2}\Big\}dt,
\end{eqnarray}
where $\rho\geq 0$.

\textbf{Problem (PI):} Seek a set of decentralized control to optimize the social cost under the worst-case disturbance
for System (\ref{eq1}) and (\ref{eq2b}), i.e.,
$\inf_{u_i\in {\cal U}_{i}}\sup_{f\in {\cal U}_{c}}J_{\rm soc},$
where $J_{\rm soc}(u,f)=\sum_{i=1}^NJ_i(u,f),$ and
$$\begin{aligned}
{\cal U}_{c} &=\Big\{u_i|\
  u_i(t)\in {\mathcal F}_t,  \mathbb{E}\int_0^\infty e^{-\rho t}\|x_i(t)\|^2dt<\infty \Big\}.
\end{aligned}
$$

\subsection{Decentralized Control Design}

Let $u_i=\check{u}_i\in {\cal U}_{c},  i=1,\cdots, N$ be fixed.
The optimal control problem with respect to drift uncertainty is as follow:
$$\hbox{\textbf{(P4)} minimize}_{f\in {\mathcal U}_{c}}\check{J}_{\rm soc}(\check{u},f),$$
where
\begin{eqnarray}\label{eq50}
\check{J}_{\rm soc}(\check{u},f)&\hspace*{-0.2cm}=\hspace*{-0.2cm}&
\frac{1}{2}\sum_{i=1}^N\mathbb{E}\int_0^{\infty}e^{-\rho t}\Big\{-\big\|x_i(t)
  -\Gamma x^{(N)}(t)\cr
  &\hspace*{-0.2cm}\hspace*{-0.2cm}&\qquad\qquad-\eta\big\|^2_{Q}
+\|f(t)\|^2_{R_2}\Big\}dt.
\end{eqnarray}

\emph{An example of scalar model}. Consider the case of uniform agents with scalar states. Let $ A=a, \eta=0, G=0,Q=1, \Gamma=\gamma$,  $R_2=r_2$, $x_i\in\mathbb{R}$ and
 $\check{u}_i=0, i=1,\cdots,N$.
By rearranging the integrand of ${J}_{\rm soc}$, we have
$$\check{J}_{\rm soc}=\frac{1}{2}\mathbb{E}\int_0^{\infty}e^{-\rho t}\Big(-x^T\hat{Q}x
+Nr_2f^2\Big)dt,$$
where $x=(x_1,\cdots,x_N)^T$, and $\hat{Q}=(\hat{q}_{ij})$ is given by
$$\hat{q}_{ii}=1+(\gamma^2-2\gamma)/N,\quad  \hat{q}_{ij}=(\gamma^2-2\gamma)/N,\ i\neq j.
$$

Introduce the Riccati equation
\begin{equation}\label{eq6}
  2(a-\frac{\rho}{2})P-\frac{1}{Nr_2}P\textbf{11}^TP-\hat{Q}=0.
\end{equation}
By observation, $P$ has the form
$p_{ij}=\left\{
\begin{array}{l}
  p \quad \hbox{if}\ i=j\\
   q \quad \hbox{if}\ i\neq j.
\end{array} \right.$
Denote $\bar{a}=a-\frac{\rho}{2}$. By solving (\ref{eq6}), we obtain the maximal solution as follows: $p= q+\frac{1}{2\bar{a}}$ and
$$q=\frac{1}{N}\left(r_2\bar{a}-\frac{1}{2\bar{a}}+\sqrt{r_2^2\bar{a}^2+r_2(\gamma^2-2\gamma-1)}\right).$$
The optimal control is given by
$$\begin{aligned}
  \check{f}=\frac{1}{N r_2}\textbf{1}^TPx
  =&\left[\bar{a}\!+\!\sqrt{\bar{a}^2\!+\!\frac{1}{r_2}(\gamma^2\!-\!2\gamma\!\!-1)}\right]x^{(N)}.
\end{aligned}$$
\rightline{$\Box$}

For general systems, we make the following assumptions:

\textbf{(A5)} Problem (P4) is uniformly convex in $f$;

\textbf{(A6)} $A+G-\frac{\rho}{2}I$ is Hurwitz.

Below are some sufficient conditions to guarantee (A5).

\begin{proposition}\label{prop4}

Let (A6) hold. (P4) is uniformly convex in $f$, i.e., (A5) holds,
if and only if one of (i)-(iv) holds.

(i) For any $f\in {\cal U}_c$, there exists $\delta>0$
such that
\begin{align*}
 &J^{\prime}_{\rm soc}(f)
=\mathbb{E}\int_0^{\infty}e^{-\rho t}\Big(-\textbf{y}^T\hat{\textbf{Q}}\textbf{y}
+{N}f^TR_2f\Big)dt\cr
&\qquad\quad\geq \delta\mathbb{E}\int_0^{\infty}e^{-\rho t}\|f\|^2dt,
\end{align*}
where $\textbf{y}\in \mathbb{R}^{nN}$ satisfies
$$   d\textbf{y} = (\check{\textbf{A}}\textbf{y}+\textbf{1}\otimes f)dt, 
   \ \textbf{y}(0)= 0.
$$

(ii)  The equation
  $$\rho{\textbf{P}}=\check{\textbf{A}}^T\textbf{P}+\textbf{P}\check{\textbf{A}}-\hat{\textbf{Q}}-
\textbf{P}(\textbf{1}\otimes I) (NR_2)^{-1}(\textbf{1}^T\otimes I)\textbf{P} $$
admits a solution such that $\check{\textbf{A}}-(\textbf{1}\otimes I) (NR_2)^{-1}(\textbf{1}^T\otimes I)\textbf{P}-\frac{\rho}{2}(I_N\otimes I)$ is Hurwitz.

(iii) The equation
\begin{align*}
\rho P=&(A+G)^TP+P(A+G)-PR_2^{-1}P+\Psi-Q
\end{align*}
admits a solution such that $A+\bar{G}-\frac{\rho}{2} I$ is Hurwitz.

(iv) The real part of any eigenvalue of $M$ is not zero,
where
$$M=\left[
\begin{array}{cc}
{A}+G-\frac{\rho}{2} I& R_2^{-1}\\
-(I-\Gamma)^TQ(I-\Gamma)&-{A}^T-G^T+\frac{\rho}{2} I
\end{array}
\right].$$
\end{proposition}
\emph{Proof.} (A5)$\Leftrightarrow$(i) follows by Lemma 1 of \cite{LZ99}.
 We now prove (ii)$\Rightarrow$(A5). If (ii) holds, then by the completion of squares technique, we can obtain
\begin{align*}
&J^{\prime}_{\rm soc}(f)\cr
=&\mathbb{E} \int_{0}^{\infty}e^{-\rho t} N\Big\|f(t)+\frac{1}{N}R_2^{-1}(\textbf{1}^T\otimes I)\textbf{P}\textbf{y}(t)\Big\|^2_{R_2} dt\geq 0.
\end{align*}
 Clearly, $ J^{\prime}_{\rm soc}(f)=0$ leads to $f(t)=-\frac{1}{N}R_2^{-1}(\textbf{1}^T\otimes I)\textbf{P}\textbf{y}(t)$, which together with $\textbf{y}(0)=0$ further implies $f(t)\equiv0$.
From \cite{HH16} we obtain that $J^{\prime}_{\rm soc}(f)$ is positive definite, which implies that (P4) is uniformly convex.
Note that $(\check{\textbf{A}}-\frac{\rho}{2}(I_N\otimes I)), I_N\otimes I)$ is stabilizable.
From (\ref{eq4b}) and (A5)$\Leftrightarrow$(i), (P4) is uniformly convex if and only if
there exists $\delta>0$
such that
 $$\int_0^{\infty}e^{-\rho t}\big(-\|(I-\Gamma)y_i\|^2_Q
+f^TR_2f\big)dt\geq \int_0^{\infty}e^{-\rho t}\delta\|f\|^2dt. $$ Following the proof of (ii)$\Rightarrow$(A5), we obtain (iii)$\Rightarrow$(A5).
 \textcolor{blue}{Since (A6) holds, it follows by \cite{SY17} that (A5)$\Rightarrow$(iii). Note that $(\textbf{1}^T\otimes I)\check{\textbf{A}}=\textbf{1}^T\otimes (A+G)$,
 $\check{\textbf{A}}(\textbf{1}\otimes I)=\textbf{1}\otimes (A+G)$, and $\frac{1}{N}(\textbf{1}^T\otimes I)\hat{\textbf{Q}}(\textbf{1}^T\otimes I)=Q-\Psi$. We have
$\frac{1}{N}(\textbf{1}^T\otimes I){\textbf{P}}(\textbf{1}^T\otimes I)=P.$
From (\ref{eq4}) and $y_1=y_2=\cdots=y^{(N)}$, we obtain  (iii)$\Leftrightarrow$(ii).}
 (iii)$\Leftrightarrow$(iv) is implied from \cite{M77}.
$\hfill \Box$

\begin{remark}
From the proof of Proposition \ref{prop4}, Assumption (A5) implies $$J^{\prime}_{soc}(f)=\mathbb{E}\int_0^{\infty}e^{-\rho t}\Big(-\textbf{y}^T\hat{\textbf{Q}}\textbf{y}
+{N}f^TR_2f\Big)dt\geq 0,$$ i.e., (P4) is convex in $f$.
\end{remark}

With some abuse of notation, later we still use $P, K, Y, Z,s, \bar{s}, \varphi,v\cdots$. But in this section $P, K, Y, Z$ are time-invariant and  $s, \bar{s}, \varphi,v$ are functions of time $t$, $t\in[0,\infty)$.
Following (\ref{eq11})-(\ref{eq11b}),
we may construct
\begin{equation*}\label{eq42}
\check{f}=-{R_2^{-1}}(Px^{(N)}+\check{s}),
\end{equation*}
where $P\in \mathbb{R}^{n\times n}$ and $\check{s}\in L^{2}_{\mathcal{F},\frac{\rho}{2}}(0, \infty; \mathbb{R}^{n})$ are determined by
\begin{align}\label{eq43e}
&(A+G-\frac{\rho}{2}I)^TP+P(A+G-\frac{\rho}{2}I)\cr
&-PR_2^{-1}P-(I-\Gamma)^TQ(I-\Gamma)=0,\\
\label{eq44}
&d\check{s}+[(A+G-R_2^{-1}P-{\rho}I)^T\check{s}+PB\check{u}^{(N)}+\bar{\eta}]dt\cr
&+\frac{1}{N}\sum_{i=1}^N\zeta_idW_i=0.
\end{align}



\begin{theorem}
  Under (A1) and (A5), Problem (P4) has a minimizer $\check{f}=-{R_2^{-1}}(Px^{(N)}+\check{s})$, where $P$ is the maximal solution of (\ref{eq43e}) and $\check{s}$ is the unique solution of (\ref{eq44}) in $L^{2}_{\mathcal{F},\frac{\rho}{2}}(0, \infty; \mathbb{R}^{n})$.
\end{theorem}
\emph{Proof.} Denote $\acute{x}_i=e^{-\frac{\rho}{2} t}x_i$, $\acute{u}_i=e^{-\frac{\rho}{2} t}u_i$ and $\acute{f}=e^{-\frac{\rho}{2}t}f$. It follows by (\ref{eq1}) and (\ref{eq50}) that
\begin{align*}
d\acute{x}_i(t) =& [(A\!-\!\frac{\rho}{2}I)\acute{x}_i(t)\!+\!B\acute{u}_i(t)+\!G\acute{x}^{(N)}(t)\!+\!\acute{f}(t)]dt\cr &+e^{-\frac{\rho}{2}t}\sigma dW_i(t),\quad 1\leq i\leq N,
\end{align*}
$$\begin{aligned}\!\!\!\!\!\!\!\!\!\!\!\!\check{J}_{\rm soc}(\check{u},f)=\frac{1}{2}\sum_{i=1}^N\mathbb{E}\int_0^{\infty}\Big\{-\big\|\acute{x}_i(t)
   -\Gamma \acute{x}^{(N)}(t)\cr-e^{-\frac{\rho}{2}t}\eta\big\|^2_{Q}
+\|\acute{f}(t)\|^2_{R_2}\Big\}dt.\end{aligned}$$
By a similar argument in the proof of Theorem \ref{thm1}, we obtain
$\check{J}_{\rm soc}(\check{u},f)= \varepsilon\Lambda^{\prime}_1+\frac{\varepsilon^2}{2}\Lambda^{\prime}_2,$
where
\begin{align*}
\Lambda^{\prime}_1\stackrel{\Delta}{=}&\sum_{i=1}^N\mathbb{E}\int_0^{\infty} \Big[\big\langle-Q\big(\acute{x}_i-(\Gamma\acute{x}^{(N)}+\eta)\big),\acute{y}_i\cr&-\Gamma \acute{y}^{(N)}\big\rangle +
\langle R_2 \acute{f},\acute{f}\rangle \Big]dt,\cr
\Lambda^{\prime}_2\stackrel{\Delta}{=}&\sum_{i=1}^N\mathbb{E}\int_0^{\infty}\Big\{-\big\|\acute{y}_i
   -\Gamma \acute{y}^{(N)}\big\|^2_{Q}
+\|\acute{f}\|^2_{R_2}\Big\}dt,
\end{align*}
and $\acute{y}_i$ satisfies
$$d\acute{y}_i=\big[(A-\frac{\rho}{2})\acute{y}_i+G\acute{y}^{(N)}+\acute{f}\big]dt, \quad \acute{y}_i(0)=0.$$
By (A5), $\Lambda^{\prime}_2\geq0$. Problem (P4) has a unique minimizer $\acute{f}=-{R_2^{-1}}\acute{p}^{(N)}$ if and only if
\begin{equation}\label{eq53b}
\left\{
\begin{aligned}
 d\acute{x}_i=&\big[(A-\frac{\rho}{2}I)\acute{x}_i+B\acute{u}_i-{R_2^{-1}}\acute{p}^{(N)}+G\acute{x}^{(N)}\big]dt\\
 &+e^{-\rho t}\sigma dW_i,\ x_i(0)={x_{i0}}, \\
d\acute{p}_i=&-\big[(A-\frac{\rho}{2}I)^T\acute{p}_i+G^T\acute{p}^{(N)}-Q\acute{x}_i+\Psi \acute{x}^{(N)}\\
&+e^{-\rho t}\bar{\eta} \big]dt+\sum_{j=1}^N\acute{\beta}_i^jdW_j,
\end{aligned}\right.
\end{equation}
admits a set of solutions $(\acute{x}_i,\acute{p}_i,i=1,\cdots,N)$ in $L^{2}_{\mathcal{F},\frac{\rho}{2}}(0, \infty; \mathbb{R}^{n})$.
 It follows from (\ref{eq53b}) that
 \begin{equation}\label{eq54b}
\left\{
\begin{aligned}
 d\acute{x}^{(N)}=&\big[(A+G-\frac{\rho}{2}I)\acute{x}^{(N)}+B\acute{u}^{(N)}\cr
 &-{R_2^{-1}}\acute{p}^{(N)}\big]dt+e^{-\rho t}\sigma dW_i,\!\ x_i(0)={x_{i0}}, \\
d\acute{p}^{(N)}=&-\big[(A+G-\frac{\rho}{2}I)^T\acute{p}^{(N)}+(\Psi -Q)\acute{x}^{(N)}\cr
&+e^{-\rho t}\bar{\eta} \big]dt+\sum_{j=1}^N\acute{\beta}_i^jdW_j.
\end{aligned}\right.
\end{equation}
 Note that $A+G-\frac{\rho}{2}I$ is Hurwitz.
By (A5) and Proposition \ref{prop4}, we obtain that (\ref{eq43e}) admits a maximal solution such that ${A}+G-\frac{\rho}{2} I-R_2^{-1}P$ is Hurwitz, which with \cite[Theorem 3.3]{SZ17} gives that
(\ref{eq44}) admits a unique solution in $L^{2}_{\mathcal{F},\frac{\rho}{2}}(0, \infty; \mathbb{R}^{n})$.
 Let $\acute{p}^{(N)}=P\acute{x}^{(N)}+\acute{s}$, where $\acute{s}=e^{-\frac{\rho}{2}t}\check{s}$. Then we have that $(\acute{x}^{(N)},\acute{p}^{(N)})$ is a solution of (\ref{eq54b}).
  By a similar argument to Theorem \ref{thm2}, the proof is completed.
\hfill$\Box$

After the worst-case drift $\check{f}$ is applied, we have the following optimal control
problem.

 \textbf{(P2$^{\prime}$)}: Minimize ${J}_{\rm soc}(u,\check{f}(u))$ over $\{u_i,1\leq i\leq N)|u_i\in {\mathcal U}_{c}\}$, where $s\in L^{2}_{\mathcal{F},\frac{\rho}{2}}(0, \infty; \mathbb{R}^{n})$,
\begin{align}
\label{eq66a}
dx_i =& [Ax_i\!+Bu_i+Gx^{(N)}\!-{R_2^{-1}}(Px^{(N)}+s)]dt\cr
&+\sigma dW_i, x_i(0)=x_{i0}, 1\leq i\leq N,\\
\label{eq66b}
ds=&-\big[(A+\bar{G})^Ts+PB{u}^{(N)}+\bar{\eta}\big]dt\cr
&+\frac{1}{N}\sum_{i=1}^N\sum_{j=1}^N(\beta_i^j-\frac{\sigma}{N})dW_j,\\
\label{eq66c}
 J_{\rm soc}(u)&=\frac{1}{2}\sum_{i=1}^N
\mathbb{E}\int_0^{\infty}e^{-\rho t}\Big\{\big\|x_i
   -\Gamma x^{(N)}-\eta\big\|^2_{Q}\cr
   &+\|u_i\|^2_{R_1}-\|Px^{(N)}+s\|^2_{R_2^{-1}}\Big\}dt.
\end{align}
\begin{lemma}\label{lem1b}
  \textcolor{blue}{Assume that A0), A5), A6) hold. 
   Then there exists ${C}_0^{\prime}>0$ such that $R_1>{C}_0^{\prime}I$ and $R_2>{C}_0^{\prime}I$, then Problem (P2$^{\prime}$) is uniformly convex.}
\end{lemma}
\emph{Proof.} Let $\textbf{z}\in \mathbb{R}^{nN}$ and $\grave{s}\in \mathbb{R}^{n}$ satisfy
\begin{align}
  \label{eq67}
   d\textbf{z} = &(\bar{\textbf{A}}\textbf{z}+\textbf{B}\textbf{u}-\textbf{1}\otimes R_2^{-1} \grave{s} )dt,\
   \textbf{z}(0)= 0,\\ \label{eq68}
   d\grave{s}=&-\big[(A+\bar{G})^T\grave{s}+\frac{1}{N}PB(\textbf{1}^T\otimes I) \textbf{u}\big]dt\cr
&+\frac{1}{N}\sum_{i=1}^N\sum_{j=1}^N\grave{\beta}_i^jdW_j.
\end{align}
  By a similar argument with \cite{LZ99}, we obtain that Problem (P2) is uniformly convex if for any $u_i\in {\mathcal U}_{c}$, there exists $\delta>0$ such that
 \begin{align}\label{eq60}
   &\mathbb{E}\int_0^{\infty}e^{-\rho t}\Big(\textbf{z}^T{\bf{\bar{Q}}}\textbf{z}
  +\textbf{u}^T\textbf{R}_1\textbf{u}-N\grave{s}^{T}R_2^{-1}\grave{{s}}\Big)dt\cr
\geq& \delta\mathbb{E} \int_0^{\infty}e^{-\rho t}\|\textbf{u}\|
^2dt.\end{align} 
Note that $A+\bar{G}-\frac{\rho}{2}I$ is Hurwitz. By \cite[Lemma 2.5]{SY17} and (\ref{eq68}),
\begin{equation}\label{eq68b}
\begin{split}
  &\mathbb{E}\int_0^{\infty}e^{-\rho t}\|\acute{s}(t)\|^2dt
\leq \frac{C_1}{N}\mathbb{E}\int_0^{\infty}e^{-\rho t} \|\textbf{u}(t)\|^2dt.
\end{split}
\end{equation}
\textcolor{blue}{Since 
${\bf{\bar{A}}}-\frac{\rho}{2}I$ is Hurwitz,
then from (\ref{eq67}) and (\ref{eq68b}), we obtain
$$\mathbb{E}\int_0^{\infty}e^{-\rho t}\|\textbf{z}\|^2dt\leq C\int_0^{\infty}e^{-\rho t}\mathbb{E}\|\textbf{u}\|^2dt.$$
Note
$  \lambda_{min}({\bf{\bar{Q}}})\geq -\lambda_{max}(\Psi+PR_2^{-1}P).
$
Thus,
 there exists $ C_0^{\prime}>0$ such that for $R_1\geq C_0^{\prime} I$ and $R_2\geq C_0^{\prime} I$, 
 (\ref{eq60}) holds.}
\hfill $\Box$

Based on the analysis in Section \ref{sec4},
we construct an auxiliary optimal control problem.

 \textbf{(P5)}: Minimize $ \bar{J}_i(u_i)$ over $u_i\in{\cal U}_i$, where
\begin{align}
 \label{eq19c1} d\grave{x}_i \!=& [A\grave{x}_i\!+\!Bu_i\!+G\bar{x}\!-{R_2^{-1}}(P\bar{x}\!+\!\bar{s})]dt+\sigma dW_i,\cr
 &\grave{x}_i(0)=x_{i0},\\
d\psi_i=&-\big[(A+\bar{G})^T \psi_i+PB u_i\big]dt
+ z_i^idW_i\cr
  \label{eq52a} \bar{J}_i(u_i)=& \frac{1}{2}\mathbb{E}\int_0^{\infty}e^{-\rho t} \Big[\grave{x}_i^TQ \grave{x}_i-2(\Psi \bar{x}+\bar{\eta})^T
\grave{x}_i\cr
&-2(P\bar{x}+\bar{s})^TR_2^{-1}P\grave{x}_i  +2v^T(\bar{G} \grave{x}_i-R_2^{-1} \psi_i)\cr
  &-2(P\bar{x}+\bar{s})^TR_2^{-1} \psi_i+u_i^TR_1u_i \Big]dt.
\end{align}
Here $\bar{s},v\in C_{\rho/2}([0,\infty),\mathbb{R}^n)$ are determined by
\begin{align*}
\dot{\bar{s}}=&-(A+\bar{G}-\rho I)^T\bar{s}-PB\bar{u}-\bar{\eta},
\cr
    \dot{v}=&-(A+\bar{G}-\rho I)^Tv+(\Psi-Q) \bar{x}+\bar{\eta}
    +PR_2^{-1}(P\bar{x}+s).
  \end{align*}
By using the 
method in \cite{YZ99, WZ12}, we can show that if (A6) holds and $Q\geq0$, (P5) admits the unique optimal control
\begin{equation}\label{eq14b}
\begin{aligned}
\hat{u}_i(t)=\ &-{R_1^{-1}}B^T(Kx_i(t)-Pl(t)+\varphi(t)),
\end{aligned}
\end{equation}
where $K\in \mathbb{R}^{n\times n}$ and $l,\varphi\in C_{\rho/2}([0,\infty),\mathbb{R}^n)$ are determined by
\begin{align*}\label{eq26c}
\rho K=&A^TK+KA-KBR_1^{-1}B^TK+Q,\\
   \dot{l}=&(A+\bar{G})l+R_2^{-1}v+R_2^{-1}(P\bar{x}+\bar{s}),\ l(0)=0,\cr
\rho \varphi=&\dot{\varphi}+\bar{A}^T\varphi+KBR_1^{-1}B^TPl+K\bar{G}\bar{x}\\
&-KR^{-1}_2\bar{s}
  -(\Psi \bar{x}+\bar{\eta})
 \!-\!PR_2^{-1}(P\bar{x}+\bar{s})+\bar{G}^Tv.
\end{align*}
By applying the control (\ref{eq14b}) into (\ref{eq66a}) combined with mean field approximations, we obtain the following equation system:
\begin{equation}\label{eq29d}
  \left\{
\begin{aligned}
\dot{\bar{x}} =& (\bar{A}\!+\!\bar{G})\bar{x}\!+\!B\!R_1^{-1}\!B^T\!(Pl\!-\!\varphi)\!-\!{R_2^{-1}}\!\bar{s}, \bar{x}(0)=\!\bar{x}_0\!,\cr
 \dot{l}=&(A+\bar{G})l+R_2^{-1}v+R_2^{-1}(P\bar{x}+\bar{s}),\ l(0)=0,\cr
\dot{\bar{s}}=&\!-\!(A+\bar{G}\! -\!\rho I)^T\bar{s}+PBR_1^{-1}B^T(K\bar{x}+\varphi)\!-\!\bar{\eta}, \cr
\dot{\varphi}\!=\!&-(\bar{A}-\!\rho I)^T\varphi-\!KBR_1^{-1}B^TPl-K\bar{G}\bar{x}\\
 &+KR^{-1}_2\bar{s}+\Psi\bar{x}+\bar{\eta}+PR_2^{-1}(P\bar{x}+\bar{s})\!-\!\bar{G}^Tv, \cr
 \dot{v}=&-(A+\bar{G} -\rho I)^Tv+(\Psi-Q)\bar{x}+\bar{\eta}\\
 &+\!PR_2^{-1}(P\bar{x}+\bar{s}).
\end{aligned}
\right.
\end{equation}

For further analysis,
we assume:

\textbf{(A7)} (\ref{eq29d}) admits a unique solution
$(\bar{x},l,s,\varphi,v)$ in $C_{\rho/2}([0,\infty), \mathbb{R}^{5n})$.

{The existence and uniqueness of a solution to (\ref{eq29d}) may be obtained by using fixed-point methods similar to those in \cite{HCM07} and \cite{WZ12}.}
We now give a sufficient condition that ensures (A7) by virtue of Riccati equations.
Using the notation in Section \ref{sec4}, we have
\begin{equation}\label{eq26e}
\left[\begin{array}{c}
 \dot{\bar{x}}\\
 \dot{{l}}\\
 \dot{s}\\
\dot{ \varphi}\\
\dot{v}
\end{array}\right]=\left[\begin{array}{cc}
M_{11}& M_{12}\\
M_{21}&M_{22}+\rho I_{3n}
\end{array}\right]
\left[\begin{array}{c}
 {\bar{x}}\\
  l\\
 {s}\\
{ \varphi}\\
{v}
\end{array}\right]+\left[\begin{array}{c}
0\\
0\\
-\bar{\eta}\\
\bar{\eta}\\
\bar{\eta}
\end{array}\right]\!\!.
\end{equation}

\begin{proposition}
If the algebraic Riccati equation
$$M_{21}+\rho Y+M_{22}Y-YM_{11}-YM_{12}Y=0$$
admits a solution $Y\in \mathbb{R}^{3n\times 2n}$ such that both $M_{11}+M_{12}Y-\frac{\rho}{2}I_{2n}$ and $-M_{22}+YM_{12}-\frac{\rho}{2}I_{3n}$ are Hurwitz, then ({A7}) holds.
  \end{proposition}
\emph{Proof.} Denote $m=[\bar{x}^T,l^T]^T$, $z=[s^T,\varphi^T,v^T]^T$. Let $z=Y m+\alpha$. By (\ref{eq26e}) and It\^{o}'s formula,
we obtain
\begin{align}\label{eq27c}
\dot{Y}=&M_{21}+M_{22}Y-YM_{11}-YM_{12}Y,\\
\label{eq27d}
\dot{\alpha}=&(M_{22}-YM_{12})\alpha+[-\bar{\eta}^T,\bar{\eta}^T,\bar{\eta}^T]^T.
\end{align}
Since (\ref{eq27c}) admits a solution such that $-M_{22}+YM_{12}-\frac{\rho}{2}I_{3n}$ is Hurwitz, then (\ref{eq27d}) has a unique solution
\begin{align*}
\alpha(t)\equiv&-\int_0^{\infty}\!\exp\big[\big(-M_{22}+YM_{12}-\frac{\rho}{2}I_{3n}\big)\tau\big]\\
&\cdot[-\bar{\eta}^T,\bar{\eta}^T,\bar{\eta}^T]^Td\tau.
\end{align*}
Applying $z=Y\bar{x}+\alpha$ into (\ref{eq29d}), we have
\begin{equation*}
\dot{m}=(M_{11}+M_{12}Y)
m+M_{12}\alpha.
\end{equation*}
Since
$M_{11}+M_{12}Y-\frac{\rho}{2}I_{2n}$ is Hurwitz, then $[\bar{x}^T, l^T]^T\in C_{\rho/2}([0,\infty), \mathbb{R}^{2n})$, and this further implies that (\ref{eq29d}) admits a unique solution in $C_{\rho/2}([0,\infty), \mathbb{R}^{5n})$.   \hfill $\Box$

\subsection{Asymptotic Optimality}

Let
\begin{equation}\label{eq70}
\begin{aligned}
\hat{u}_i(t)=\ &-{R_1^{-1}}B^T(Kx_i(t)-Pl(t)+\varphi(t)),
\end{aligned}
\end{equation}
where $ l$ and $\varphi$ are determined by (\ref{eq29d}).
After the control $\hat{u}_i$ is applied, the closed-loop dynamics can be written as
\begin{equation}\label{eq31c}
   \begin{split}
 d\hat{x}_i =[\bar{A}\hat{x}_i+\bar{G}\hat x^{(N)}+BR_1^{-1}B^T(Pl-\varphi)\cr
 -{R_2^{-1}}\hat{s}]dt+\sigma dW_i. 
   \end{split}
\end{equation}

For further analysis,
we assume

\textbf{(A8)} The equation
\begin{equation}\label{eq70a}
   \tilde{P}(\bar{A}+\bar{G})+({A}+\bar{G})^T\tilde{P}-\tilde{P}R_2^{-1}\tilde{P}+PBR_1^{-1}B^TK=0
   \end{equation}
admits a solution $\tilde{P}$ such that
$\bar{A}+\bar{G}-\frac{\rho}{2}I-R_2^{-1}\tilde{P}$ and ${A}+\bar{G}-\frac{\rho}{2}I-R_2^{-1}\tilde{P}$ are Hurwitz, where $\bar{A}=A-BR_1^{-1}B^TK$ and $\bar{G}=G -R_2^{-1}P$.

\begin{theorem}\label{thm3.2}
Assume (i) (A0)-(A1), (A5)-(A8) hold, (ii) $A-\frac{\rho}{2}I$ is Hurwitz (iii) (P2$^{\prime}$) is convex. For Problem (PI), the set of control laws
$\hat{u}=(\hat{u}_1,\cdots,\hat{u}_N)$ given by (\ref{eq24b}) has asymptotic robust social optimality, i.e.,
\begin{align*}
&\left|\frac{1}{N} \sup_{f\in {\mathcal U}_c}J_{\rm soc}(\hat{u},f)-\frac{1}{N}\inf_{u_i\in{\mathcal U}_c }
\sup_{f\in {\mathcal U}_c}J_{\rm soc}({u},f) \right|
=O(\frac{1}{\sqrt{N}}).
\end{align*}
\end{theorem}

\emph{Proof.} See Appendix \ref{app b}.
\hfill$\Box$

\section{Numerical Example}

\begin{figure}[H]
	\centering
{	\includegraphics[width=0.8\linewidth]{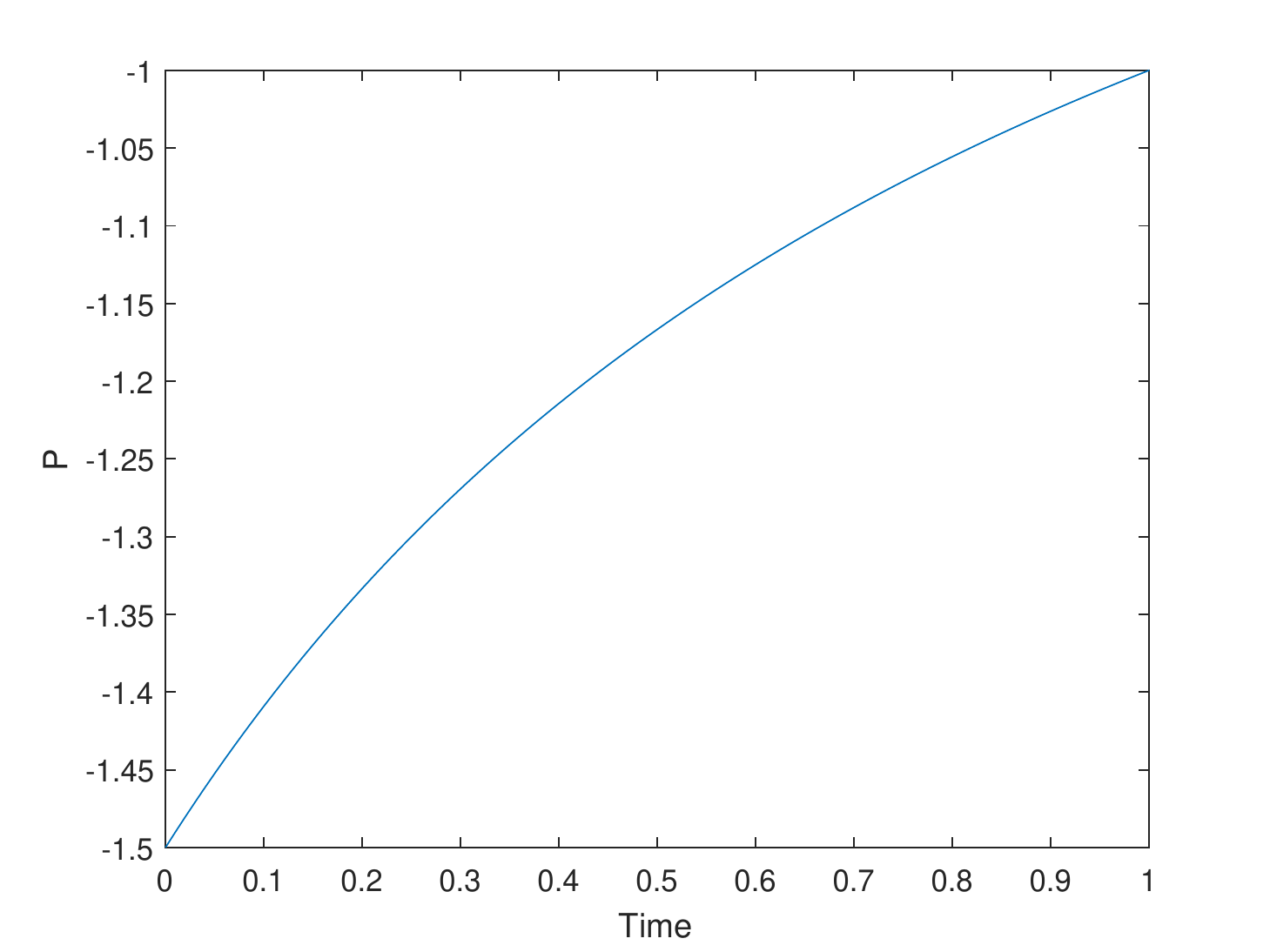}}
	\caption{The curve of $P(t)$}
	\label{fig:P}
\end{figure}

\begin{figure}[H]
	\centering
	\hspace*{-0.4cm}\includegraphics[width=1.0\linewidth]{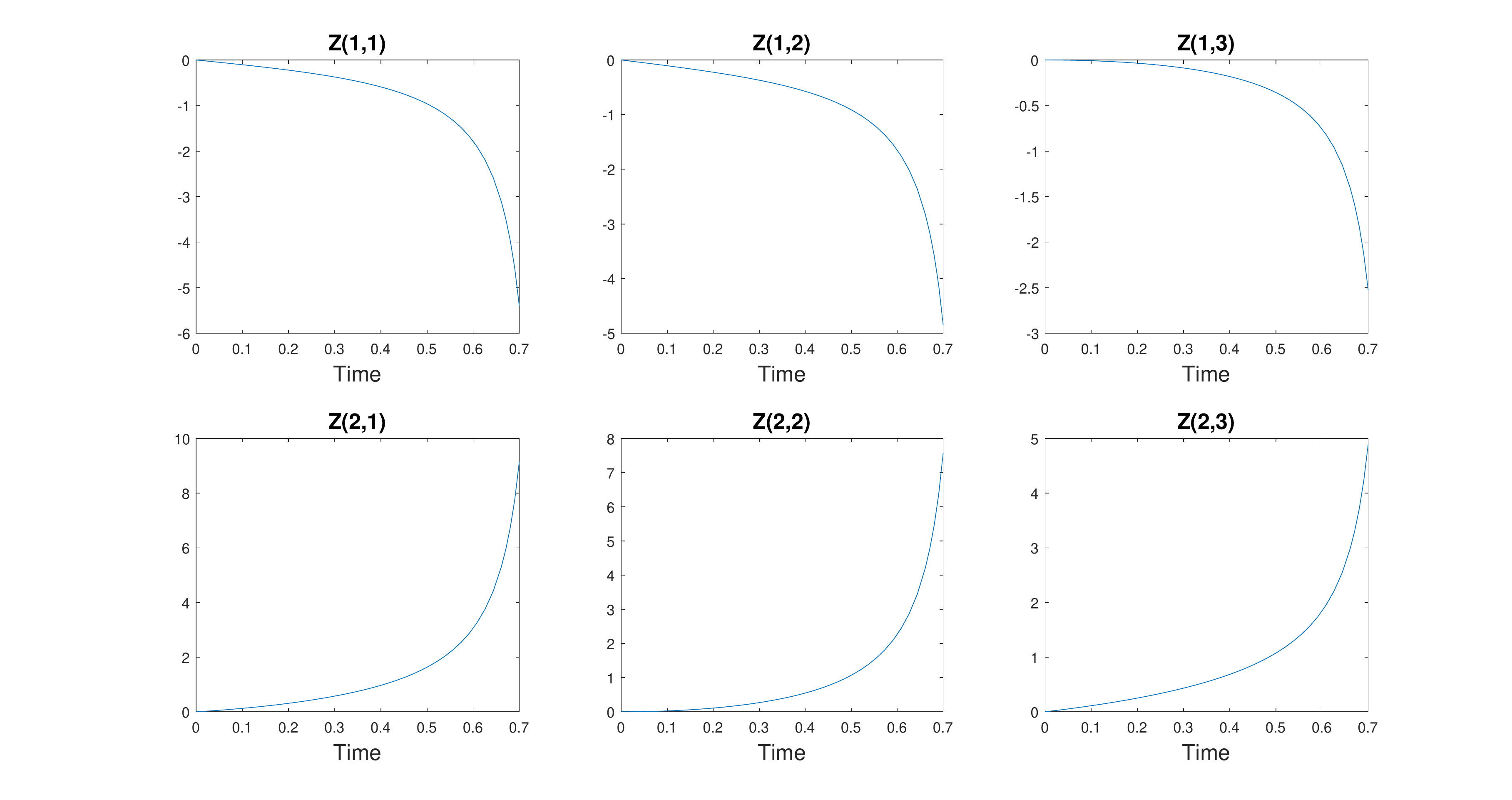}
	\caption{The curves of all entries of $Z \in \mathbb{R}^{2\times3}$ when $t \in [0,0.7]$}
	\label{fig:untitled}
\end{figure}

\begin{figure}[H]
	\centering
	\includegraphics[width=0.8\linewidth]{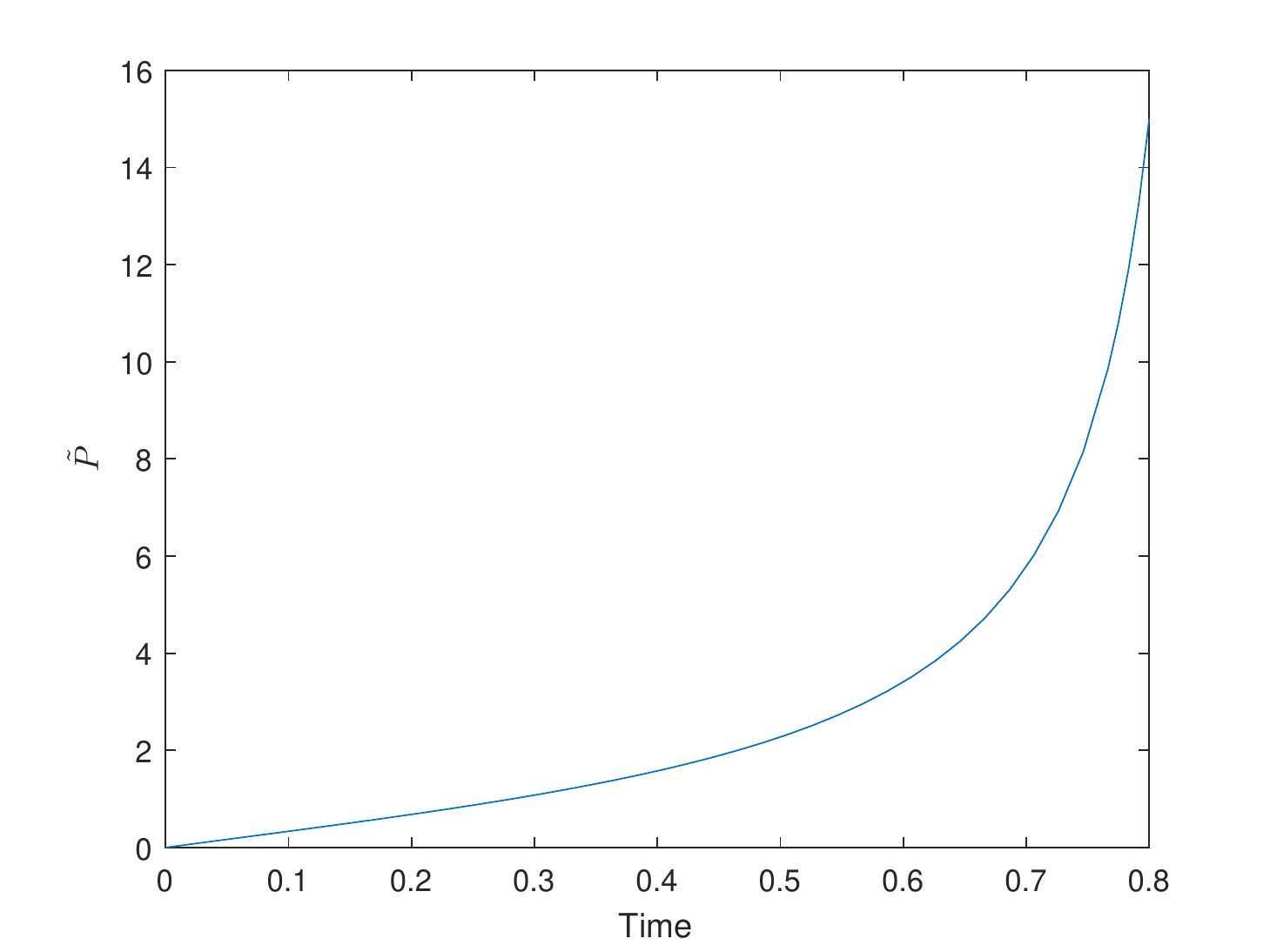}
	\caption{The curve of $\tilde{P}(t)$}
	\label{fig:tildeP}
\end{figure}

We now give a numerical example for Problem \textbf{(PF)} to verify the result. Take the parameters $A=B=R_1=R_2=Q=H=1$, $G=-1.5$, $\Gamma=0.5$, $\eta=0$, and $T=1$.
By solving (\ref{eq11}), we can obtain that $P(t)=-\frac{1}{\left(t+1\right)}-\frac{1}{2}$, which is shown in Fig. \ref{fig:P}.
 By Proposition \ref{prop1b}, (A2$^{\prime}$) holds.
For (\ref{eq38a}) in Proposition \ref{prop4a}, the curves of all entries of the solution $Z$ are given in Fig. \ref{fig:untitled}.
It can be seen that when $t\in [0,0.7]$, (\ref{eq38a}) admits a solution. By Matlab computation, the solution blows up at $t=0.758276$. From Proposition \ref{prop4a},
when $t\in [0,0.7]$, (A3) holds. The curve of $\tilde{P}$ is shown in Fig. \ref{fig:tildeP}. It can seen that  the Riccati equation in (A4) adimits a solution $\tilde{P}$  when $t \in [0,0.8]$. As a conclusion, when $t\in [0,0.7]$, (A0)-(A1), (A2$^{\prime}$), (A3)-(A4) hold. By Theorem \ref{thm4}, Problem (PF) admits a set of control laws which has asymptotic robust social optimality.

\section{Concluding Remarks}
This paper considered a class of mean field LQG social optimum problem with global drift uncertainty. 
Based on the soft control approach, a set of decentralized strategies is designed by optimizing the worst-case cost subject to consistent requirements in mean field approximations.
Such set of strategies is further shown to be robust social optimal by perturbation analysis.

For further work, it is of interest to consider mean field team optimization with volatility-uncertain common noise. Due to
common noise and volatility uncertainty, all states
of agents are coupled via some high-dimensional FBSDE systems. 
Another interesting topic is the mean field Stackelberg game with a leader and many followers \cite{WZ14, XSZ2018}. The team problem with hierarchical structure is worthwhile to study further.

\appendices
\section{Proof of Theorem \ref{thm4}}\label{app a}
\def\theequation{A.\arabic{equation}}
\setcounter{equation}{0}

\emph{Proof.} Note that we only need to optimize the social cost under worst-case disturbance $J_{\rm soc}^{\rm wo}({u})$.  By Theorem \ref{thm2}, we can restrict to consider Problem (P2) instead of (PF).
From Lemma \ref{lem4}, one obtains that for (P2),
\begin{equation}\label{eq37}
  \mathbb{E}\int_0^T(\|\hat{x}_i\|^2+\|\hat{u}_i\|^2)dt<C.
  \end{equation}
It suffices to consider all $u_i\in {\mathcal U}_c^{\rm F}$ such that $J_{\rm soc}^{\rm wo}({u})\leq J_{\rm soc}^{\rm wo}(\hat{u})\leq NC_0. $
By Lemma \ref{lem2},
\begin{equation}\label{eq38}
\mathbb{E}\int_0^T\|{u}_i\|^2dt<C,\ i=1,\cdots,N,
\end{equation}
which implies
\begin{equation}\label{eq39}
  \mathbb{E}\int_0^T\|{u}^{(N)}\|^2dt<C.
  \end{equation}
  By (\ref{eq12b}) and \cite[Chaper 7]{YZ99}, we have
\begin{equation*}\label{eq40a}
\mathbb{E}\int_0^T\|s\|^2dt\leq  C_1\mathbb{E}\int_0^T\|u^{(N)}\|^2dt<C.
\end{equation*}
From (\ref{eq12a}),
$$dx^{(N)}\!=\!\big[(A\!+\!\bar{G})x^{(N)}\!+\!Bu^{(N)}\!-\!R_2^{-1}s\big]dt\!+\!\frac{1}{N}\!\sum_{i=1}^N\sigma dW_i,$$
which together with (\ref{eq39}) implies $\mathbb{E}\int_0^T\|{x}^{(N)}\|^2dt\!<C\!\!.$
This with (\ref{eq38}) leads to
\begin{equation}\label{eq41}
\mathbb{E}\int_0^T(\|{x}_i\|^2+\|{u}_i\|^2+\|s\|^2)dt<C.
\end{equation}
Let $\tilde{x}_i=x_i-\hat{x}_i$,   $\tilde{u}_i=u_i-\hat{u}_i$, $i=1,\cdots,N$, $\tilde{x}^{(N)}=\frac{1}{N}\sum_{j=1}^N \tilde{x}_j$ and
$\tilde{s}=s-\hat{s}$. Then by (\ref{eq12a}),
\begin{align}\label{eq32}
  d\tilde{x}_i=&(A\tilde{x}_i+\bar{G}\tilde{x}^{(N)}+B\tilde{u}_i-R_2^{-1}\tilde{s})dt, \  \tilde{x}_i(0)=0.\\
 \label{eq32b} d\tilde{s}=&-\big[(A+\bar{G})^T\tilde{s}+PB\tilde{u}^{(N)}\big]dt\cr
&+\frac{1}{N}\sum_{i=1}^N\sum_{j=1}^N\tilde{\beta}_i^jdW_j
,\ \tilde{s}(T)=0.
\end{align}
By (\ref{eq37}) and (\ref{eq41}),
\begin{equation}\label{eq43}
\mathbb{E}\int_0^T(\|\tilde{x}_i\|^2+\|\tilde{u}_i\|^2+\|\tilde{s}\|^2)dt<C.
\end{equation}
From (\ref{eq12b}), we have
\begin{equation}\label{eq54}
\begin{aligned}
    J^{\rm F}_{\rm soc}(u)= &\frac{1}{2} \sum_{i=1}^N\mathbb{E}\int_0^T\Big[\big\|\hat{x}_i-\Gamma \hat{x}^{(N)}-\eta+\tilde{x}_i\cr
    &-\Gamma\tilde{x}^{(N)}\big\|^2_Q+\big\|\hat{u}_i+\tilde{u}_i\big\|^2_{R_1}\cr
    &-\big\|P(\hat{x}^{(N)}+\tilde{x}^{(N)})+\hat{s}+\tilde{s}\big\|^2_{R_2^{-1}}\Big]dt\cr
    &+\frac{1}{2}\mathbb{E}\big\|\hat{x}_i(T)+\tilde{x}_i(T)\big\|^2_H
\\    =&\sum_{i=1}^N(J_i^{\rm F}(\hat{u})+\tilde{J}_i^{\rm F}(\tilde{u})+I_i),
\end{aligned}
\end{equation}
where
\begin{align*}
  \tilde{J}_i^{\rm F}(\tilde{u})\stackrel{\Delta}{=}&\frac{1}{2}\mathbb{E}\int_0^T\big[\|\tilde{x}_i-\Gamma \tilde{x}^{(N)}\|^2_Q+\|\tilde{u}_i\|^2_{R_1}\\
  &-\|P\tilde{x}^{(N)}\|^2_{R_2}\big]dt +\frac{1}{2}\mathbb{E}\big\|\tilde{x}_i(T)\big\|^2_H,\\
  I_i=&\mathbb{E}\int_0^T\Big[\big(\hat{x}_i-\Gamma \hat{x}^{(N)}-\eta\big)^TQ\big(\tilde{x}_i-\Gamma \tilde{x}^{(N)}\big)\\
    &+\hat{u}_i^TR_1\tilde{u}_i-\big(P(\hat{x}^{(N)}+\hat{s})\big)^T{R_2^{-1}}(P\tilde{x}^{(N)}\\
    &+\tilde{s})\Big]dt +\mathbb{E}[\hat{x}_i^T(T)H\tilde{x}_i(T)].
\end{align*}
By Lemma \ref{lem1a}, Problem (P2) is uniformly convex for $N\geq N_0$, which with Proposition \ref{prop1} gives $\tilde{J}_i(\tilde{u})\geq 0$. We now prove $\frac{1}{N}\sum_{i=1}^N I_i=O(\frac{1}{\sqrt{N}})$. By straightforward computation,
\begin{align}\label{eq40}
  &\sum_{i=1}^N I_i=\sum_{i=1}^N \mathbb{E}\int_0^T \Big\{\tilde{x}_i^T\big[Q(\hat x_i-\Gamma \hat{x}^{(N)}-\eta)\cr
  &-\Gamma^TQ((I-\Gamma)\hat{x}^{(N)}-\eta)\big]+\hat{u}_i^TR_1\tilde{u}_i\Big\}dt\cr
  &-N\mathbb{E}\int_0^T\big(P\hat{x}^{(N)}+\hat{s}\big)^T{R_2^{-1}}
  (P\tilde{x}^{(N)}+\tilde{s})dt\cr
  &+ \sum_{i=1}^N\mathbb{E}[\hat{x}_i^T(T)H\tilde{x}_i(T)]\cr
  =&\sum_{i=1}^N\mathbb{E}\int_0^T \Big\{\tilde{x}_i^T(Q\hat x_i-\Psi\bar{x}-\eta)+\hat{u}_i^TR_1\tilde{u}_i\cr
  &-\big(P\bar{x}+\hat{s}\big)^T{R_2^{-1}}P\tilde{x}_i\Big\}dt\cr
  &+\sum_{i=1}^N\mathbb{E}\!\int_0^T\!\!\xi^T\big[(\Psi \!-\!PR_2^{-1}P)\tilde{x}_i\!-\!PR_2^{-1}\tilde{s}\big]dt\cr
 & -N\mathbb{E}\!\!\int_0^T\!\!(P\bar{x}+\hat{s})^TR_2^{-1}\tilde{s}dt+ \sum_{i=1}^N\mathbb{E}\|\hat{x}_i^T(T)\|_H^2,
\end{align}
where $\xi=\hat{x}^{(N)}-\bar{x}$.
By (\ref{eq29}) and (\ref{eq24c}),
\begin{align}
   d\hat{k}_i=&\big\{-A^T\hat{k}_i-Q\hat{x}_i+(\Psi\bar{x}+\bar{\eta})+PR_2^{-1}(P\bar{x}+\bar{s})\cr
   &-\!\bar{G}^Tv\!+\!KG(\hat{x}^{(N)}-\bar{x})-KR_2^{-1}(\hat{s}-\bar{s})\big\}dt\cr
   &+K\sigma dW_i, \ \hat{k}_i(T)=H\hat{x}_i(T).
\end{align}
By (\ref{eq32}) and It\^{o}'s formula,
$$  \begin{aligned}
&\sum_{i=1}^N\mathbb{E}[\hat{x}_i^T(T)H\tilde{x}_i(T)]\cr
=&\sum_{i=1}^N\mathbb{E}\int_{0}^{T}\Big\{-\big[ Q\hat{x}_i-(\Psi \bar{x}+\bar{\eta})-PR_2^{-1}(P\bar{x}+\bar{s})\cr
&+\bar{G}^Tv-KG(\hat{x}^{(N)}-\bar{x})+KR_2^{-1}(\hat{s}-\bar{s})\big]^T \tilde{x}_i \cr
 &+\hat{k}_i^T(\bar{G}\tilde{x}^{(N)}+B\tilde{u}_i-R_2^{-1}\tilde{s})\Big\}dt,
\end{aligned}$$
and
$$\begin{aligned}
  0=&\sum_{i=1}^N\mathbb{E}[l^T(T)\tilde{s}(T)-l^T(0)\tilde{s}(0)]\\
  =&N\mathbb{E}\int_0^T[R_2^{-1}v+R_2^{-1}(P\bar{x}+\bar{s})]^T\tilde{s}dt\cr
  &-N\mathbb{E}\int_0^T(l^TPB\tilde{u}^{(N)})dt.
\end{aligned}$$
The above two equations lead to
$$  \begin{aligned}
&\sum_{i=1}^N\mathbb{E}[\hat{x}_i^T(T)H\tilde{x}_i(T)]\cr
=&\sum_{i=1}^N\mathbb{E}\int_{0}^{T}\Big\{-\big[ Q\hat{x}_i-(\Psi \bar{x}+\bar{\eta})-PR_2^{-1}(P\bar{x}+\bar{s})\cr
&-KG(\hat{x}^{(N)}-\bar{x})+KR_2^{-1}(\hat{s}-\bar{s})\big]^T \tilde{x}_i \cr
&- \hat{u}_i^TR_1\tilde{u}_i+(\hat{k}^{(N)}-v)^T\bar{G}\tilde{x}_i\cr
&-(\hat{k}^{(N)}-v)R_2^{-1}\tilde{s}+(P\bar{x}+\bar{s})^TR_2^{-1}\tilde{s}\Big\}dt.
\end{aligned}$$
From this and (\ref{eq40}),
 $$ \begin{aligned}
 \sum_{i=1}^N I_i=&\sum_{i=1}^N\mathbb{E}\int_{0}^{T}\big[\xi^T(\Psi -PR_2^{-1}P+KG)\tilde{x}_i\cr
 &+(\hat{k}^{(N)}-v)^T(\bar{G}\tilde{x}_i-R_2^{-1}\tilde{s})\cr
 &+ ((P+K)\tilde{x}_i+\tilde{s})R_2^{-1}(\bar{s}-\hat{s})\big]dt.
 \end{aligned}$$
By Lemmas \ref{lem1}, \ref{lem5} and Schwarz inequality, we obtain
$$\frac{1}{N}\sum_{i=1}^N I_i=O(\frac{1}{\sqrt{N}}).$$
From this with (\ref{eq54}), the theorem follows. $\hfill{\Box}$

\section{Proof of Theorem \ref{thm3.2}}\label{app b}
\def\theequation{B.\arabic{equation}}
\setcounter{equation}{0}

To prove Theorem \ref{thm3.2}, we need three lemmas.

\begin{lemma}\label{lem6}
Assume that (A0)-(A1), (A5)-(A8) hold. For Problem (PI), we have 
\begin{equation}\label{eq52}
 \mathbb{E}\int_0^{\infty}e^{-\rho t}\big(\|\hat{x}^{(N)}-\bar{x}\|^2+\|\hat{s}-\bar{s}\|^2\big)dt=O(\frac{1}{N}).
  \end{equation}
\end{lemma}
\emph{Proof.}  
By a similar argument to (\ref{eq45a})-(\ref{eq45b}), we obtain
\begin{align*}
  d\xi=&(\bar{A}+\bar{G})\xi dt-{R_2^{-1}}\chi dt+\frac{1}{N}\sum_{i=1}^N\sigma dW_i,\\
  & \textstyle{\xi(0)=\frac{1}{N}\sum_{i=1}^Nx_{i0}-\bar{x}_0},\\
   d\chi=&-\big[(A+\bar{G})^T\chi+PBR_1^{-1}B^TK\xi\big]dt\\
&+\frac{1}{N}\sum_{i=1}^N\sum_{j=1}^N(\beta_i^j-\frac{\sigma}{N})dW_j,
\end{align*}
where $\xi{=}\hat{x}^{(N)}-\bar x$ and $\chi{=}\hat{s}-\bar{s}$.
By It\^{o}'s formula and (A8), we have $\chi=\tilde{P}\xi+\psi$,
where $\tilde{P}$ is given by (\ref{eq70a}).
 Denote $\Upsilon=\bar{A}+\bar{G}-R_2^{-1}\tilde{P}$. Then
$$\begin{aligned}
\xi(t)=&e^{\Upsilon t}\xi(0)+\frac{1}{N}\int_0^te^{\Upsilon(t-\mu)}\sum_{i=1}^N
\sigma dW_i(\mu).
\end{aligned}$$
This with (A8) gives $\mathbb{E}\int_0^{\infty}e^{-\rho t}\|\xi(t)\|^2 dt=O({1}/{N})$.
\hfill $\Box$

\begin{lemma}\label{thm2b}
	Assume that (A0)-(A1), (A5)-(A8) hold. For Problem (PI) and any $N$,
	\begin{equation}\label{eq13c}
	\max_{1\leq i\leq N} \mathbb{E}\int_0^{\infty} e^{-\rho t} \left(\|\hat{x}_i\|^2+\|\hat{u}_i\|^2\right)dt<\infty.
	\end{equation}
\end{lemma}

\emph{Proof.} 
 By (A7) and Lemma \ref{lem6}, we obtain that
$$\mathbb{E}\int_0^{\infty} e^{-\rho t} (\|\hat{x}^{(N)}(t)\|^2+\|\hat{s}(t)\|^2)dt<\infty.$$
Note that $\bar{A}-\frac{\rho}{2}I$ is Hurwitz. By Schwarz's inequality,
\begin{align*}
&\mathbb{E}\int_0^{\infty}e^{-\rho t}\|\hat{x}_i(t)\|^2dt\cr
\leq\ &C\!+\!3\mathbb{E}\!\int_0^{\infty}\!\!e^{-\rho \mu}\|g(\mu)\|^2\!\int_\mu^{\infty}\!\!t\big \|e^{(\bar{A}-\frac{\rho}{2}I)(t-\mu)}\big\|^2dtd\mu\cr
&+3 C\mathbb{E}\!\int_0^{\infty}\!\!e^{-\rho \mu}\|\sigma(\mu)\|^2\!\int_\mu^{\infty}\!\!\big \|e^{(\bar{A}-\frac{\rho}{2}I)(t-\mu)}\sigma\big\|^2dtd\mu\cr
\leq\ &C+3C\mathbb{E}\int_0^{\infty}e^{-\rho \mu}\|g(\mu)\|^2d\mu\\
&+3 C\mathbb{E}\int_0^{\infty}e^{-\rho \mu}\|\sigma(\mu)\|^2d\mu\leq C_1.
\end{align*}
This with (A7) completes the proof.   \hfill $\Box$

\begin{lemma}\label{lem8}
Assume $A-\frac{\rho}{2}I$ is Hurwitz. Then
$$ \mathbb{E} \int_0^{\infty}e^{-\rho t}\|\hat{k}^{(N)}-v\|^2dt\leq O(1/N).$$ 
\end{lemma}
\emph{Proof.}
By (\ref{eq29d}) and some elementary computations, we obtain
$d\vartheta(t)=-(A-\frac{\rho}{2}I)^T\vartheta(t) dt,$
where $\vartheta=v-K\bar{x}-\varphi$. This leads to
$\vartheta(t)=e^{-(A-\frac{\rho}{2}I)t}\vartheta(0)$. Since $A-\frac{\rho}{2}I$ is Hurwitz, and $\vartheta\in C_{\rho/2}([0,\infty), \mathbb{R}^n)$, then we have $\vartheta(t)\equiv 0$, which implies $v=K\bar{x}+\varphi$.
By Lemma \ref{lem6},
$\int_0^{\infty}e^{-\rho t}\mathbb{E}\|\hat{k}^{(N)}-v\|^2dt\leq O(1/N).
$
This completes the proof.  \hfill$\Box$

\emph{Proof of Theorem \ref{thm3.2}.} As in the proof of Theorem \ref{thm4}, we restrict to Problem (P2$^{\prime}$).
It suffices to consider all $u_i\in {\mathcal U}_c$ such that
$\sup_{f\in  {\mathcal U}_c}J_{\rm soc}({u},f)\leq\sup_{f\in  {\mathcal U}_c}J_{\rm soc}(\hat{u},f)\leq NC_0.
$            Taking $f=0$, we have
\begin{equation}\label{eq38c}
\mathbb{E}\int_0^{\infty}e^{-\rho t}\|{u}_i\|^2dt<C.
\end{equation}
  By (\ref{eq66b}) and \cite{SY17}, we have
\begin{equation*}\label{eq40c}
\mathbb{E}\int_0^{\infty}e^{-\rho t}\|s\|^2dt\leq  C_1\mathbb{E}\int_0^{\infty}e^{-\rho t}\|u^{(N)}\|^2dt<C.
\end{equation*}
Noticing $A+\bar{G}-\frac{\rho}{2} I$ is Hurwitz, one can obtain $\mathbb{E}\int_0^{\infty}e^{-\rho t}\|{x}^{(N)}\|^2dt<C$
which with (\ref{eq38c}) implies
           \begin{equation}\label{eq38b}
\mathbb{E}\int_0^{\infty}e^{-\rho t}\big(\|{x}_i\|^2+\|{u}_i\|^2+\|s\|^2\big)dt\leq C.
\end{equation}
From this and (\ref{eq13c}),
\begin{equation}\label{eq43c}
\mathbb{E}\int_0^{\infty}(\|\tilde{x}_i\|^2+\|\tilde{u}_i\|^2+\|\tilde{s}\|^2)dt<C.
\end{equation}
We have
$    J_{\rm soc}(u)=\sum_{i=1}^N(J_i(\hat{u})+\tilde{J}_i(\tilde{u})+\mathcal{I}_i),
$
where
\begin{align*}
  \tilde{J}_i(\tilde{u})\stackrel{\Delta}{=}\frac{1}{2}\mathbb{E}&\int_0^{\infty}e^{-\rho t}\big[\|\tilde{x}_i-\Gamma \tilde{x}^{(N)}\|^2_Q\cr
  &+\|\tilde{u}_i\|^2_{R_1}-\|P\tilde{x}^{(N)}\|^2_{R_2}\big]dt,\\
  \mathcal{I}_i\stackrel{\Delta}{=}\mathbb{E}&\int_0^{\infty}e^{-\rho t}\Big[\big(\hat{x}_i-\Gamma \hat{x}^{(N)}-\eta\big)^TQ\big(\tilde{x}_i\cr
  &-\Gamma \tilde{x}^{(N)}\big)+\hat{u}_i^TR_1\tilde{u}_i\cr
  &-\big(P(\hat{x}^{(N)}+s)\big)^T{R_2^{-1}}(P\tilde{x}^{(N)}+\tilde{s})\Big]dt.
\end{align*}
 From Lemma \ref{lem1b} and Proposition \ref{prop4},  $\tilde{J}_i(\tilde{u})\geq 0$ for $N\geq N_0$. 
By making use of It\^{o}'s formula and straightforward computations,
 $$ \begin{aligned}
 \sum_{i=1}^N \mathcal{I}_i=&
 \sum_{i=1}^N\mathbb{E}\int_0^{\infty}e^{-\rho t} \Big\{\tilde{x}_i^T\big[Q\hat x_i-\Psi \bar{x}-\bar{\eta}\big]\cr
  & +\hat{u}_i^TR_1\tilde{u}_i-\big(P\bar{x}+\hat{s}\big)^T{R_2^{-1}}P\tilde{x}_i\Big\}dt\cr
 &+\sum_{i=1}^N\mathbb{E}\int_0^{\infty}e^{-\rho t}\xi^T[(\Psi-PR_2^{-1}P)\tilde{x}_i-PR_2^{-1}\tilde{s}]dt\cr
  &-N\mathbb{E}\int_0^{\infty}e^{-\rho t}(P\bar{x}+\hat{s})^TR_2^{-1}\tilde{s}dt,\cr
=&\sum_{i=1}^N\mathbb{E}\int_{0}^{\infty}e^{-\rho t}\big[\xi^T(\Psi -PR_2^{-1}P+KG)\tilde{x}_i\\
 &+(\hat{k}^{(N)}-v)^T(\bar{G}\tilde{x}_i-R_2^{-1}\tilde{s})\cr
 &+ ((P+K)\tilde{x}_i+\tilde{s})R_2^{-1}(\bar{s}-\hat{s})\big]dt.
 \end{aligned}$$
From (\ref{lem6}) and (\ref{eq43c}),
we obtain
$$\frac{1}{N}\sum_{i=1}^N \mathcal{I}_i=O(\frac{1}{\sqrt{N}}).$$  \hfill{$\Box$}

\vspace{-10ex}

\begin{IEEEbiography}[{\includegraphics[width=1in,height=1.25in, clip, keepaspectratio]{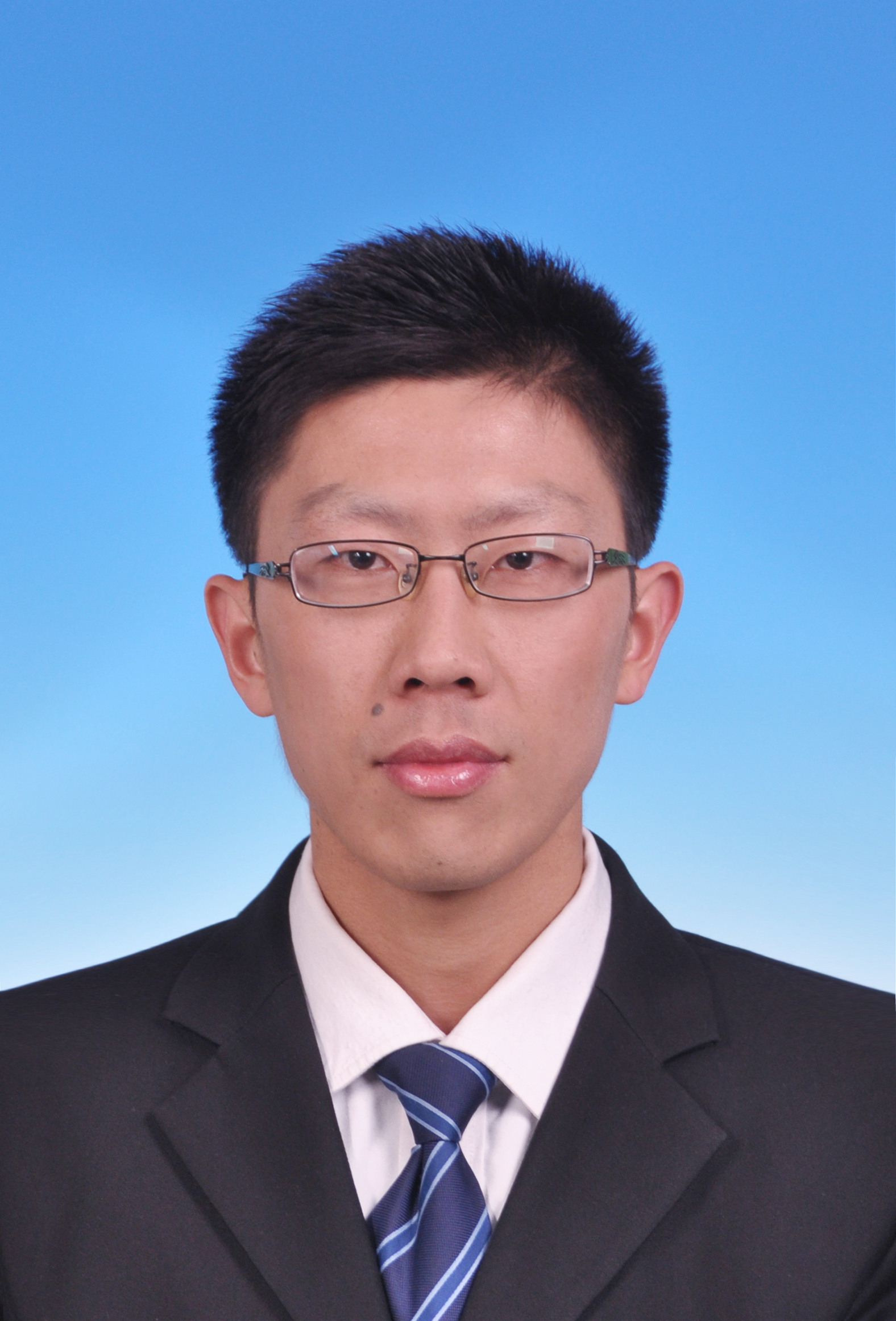}}]
{Bingchang Wang}
 received the M.Sc. degree in Mathematics from
Central South University, Changsha, China, in 2008, and
the Ph.D. degree in System Theory from Academy of
Mathematics and Systems Science, Chinese Academy of
Sciences, Beijing, China, in 2011. From September 2011
to August 2012, he was with Department of Electrical
and Computer Engineering, University of Alberta, Canada,
as a Postdoctoral Fellow. From September 2012 to September 2013, he was with School of Electrical Engineering and Computer
Science, University of Newcastle, Australia, as a Research
Academic.
From October 2013, he has
been with School of Control Science and Engineering, Shandong University, China, as an associate Professor. He held visiting
appointments as a Research Associate with Carleton University, Canada, from November 2014 to May 2015, and with the Hong Kong Polytechnic University from November 2016 to January
 2017. 
His current research interests include mean field games, stochastic control, multiagent
systems and event based control. He received the IEEE CSS Beijing Chapter Young Author Prize in 2018.
\end{IEEEbiography}

\vspace{-10ex}

\begin{IEEEbiography}[{\includegraphics[width=1in,height=1.25in, clip, keepaspectratio]{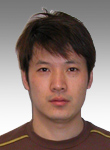}}]
{Jianhui
Huang}
Jianhui Huang received his B.S. degree in operational
research and control theory, M.Sc. degree in
probability theory and mathematical statistics from
the School of Mathematics and System Sciences,
Shandong University, Jinan, China, respectively, in
1998 and 2001. He received his Ph.D. degree in
mathematical finance from the Department of Mathematical
and Statistical Sciences, University of Alberta,
Edmonton, Canada, in 2007. Since 2007, he
has been with The Hong Kong Polytechnic University,
where he is currently an assistant professor
at the Department of Applied Mathematics. His current research interests
include stochastic control and optimization, stochastic dynamic games, large population
systems and their applications. He has published in journals such as
SIAM Journal on Control and Optimization, IEEE Transactions on Automatic
Control and Automatica.

\end{IEEEbiography}
\vspace{-10ex}

\begin{IEEEbiography}[{\includegraphics[width=1in,height=1.25in, clip, keepaspectratio]{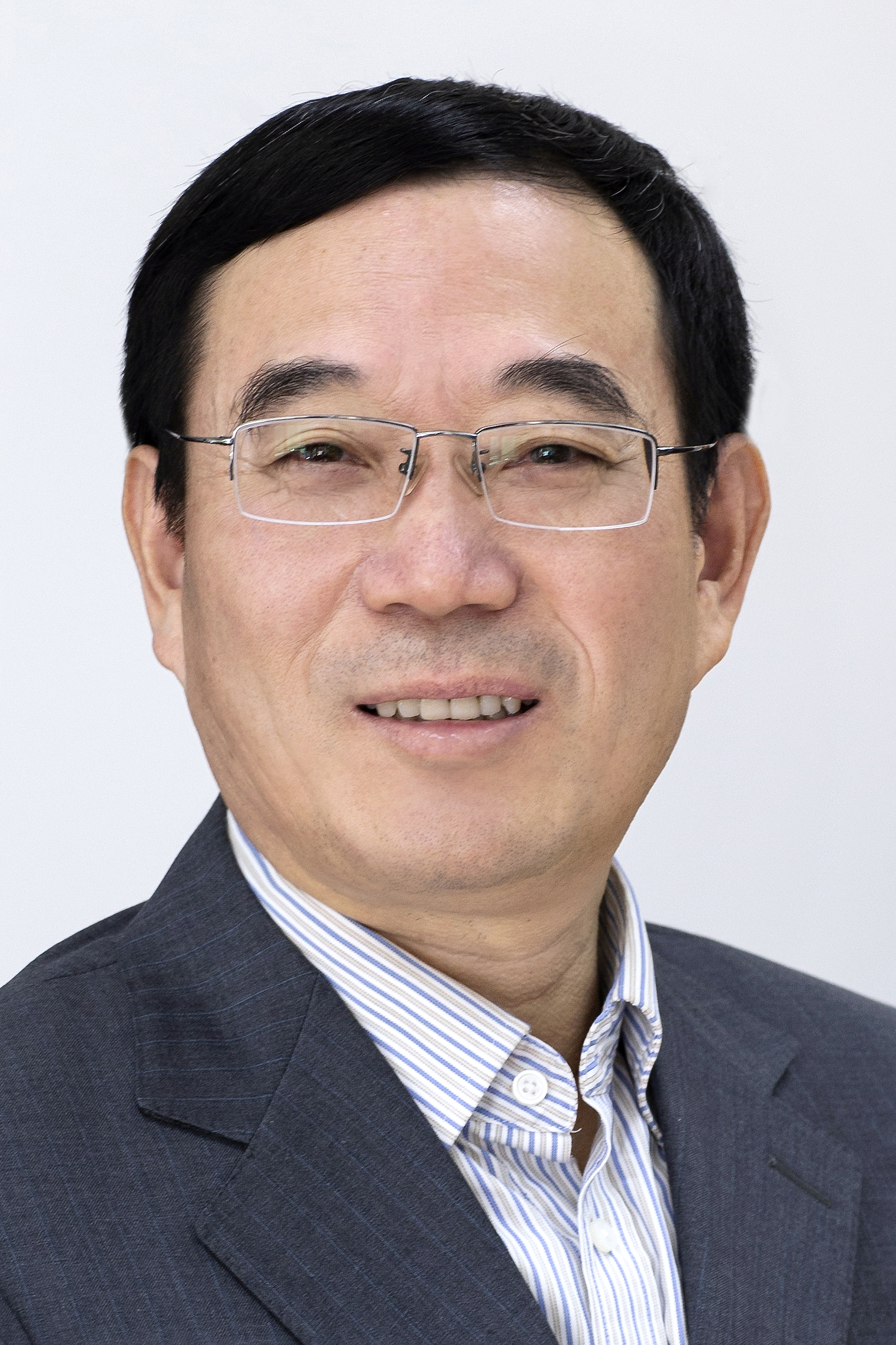}}]
{Ji-Feng Zhang}
(M'92-SM'97-F'14)
received the B.S. degree in mathematics from
Shandong University, China,
in 1985 and the Ph.D. degree from
the Institute of Systems Science (ISS),
Chinese
Academy of Sciences (CAS), China, in 1991.
Since 1985, he has been with the ISS, CAS, and now is the Director
 of ISS.
He is a
Deputy Editor-in-Chief of {\em Science China Information Sciences},
and  was the Editor-in-Chief of {\em Journal of Systems Science and
 Mathematical
Sciences}, the founding Editor-in-Chief of {\em All About Systems and Control},
a
 Managing Editor of {\em Journal of Systems Science and
 Complexity},
Deputy Editor-in-Chief of {\em Acta Automatica Sinica}, {\em Control Theory and Applications}
and {\em Systems Engineering --- Theory and Practice}, Associate Editor of
several
other journals, including {\em IEEE Trans. on
Automatic Control},
{\em SIAM Journal on Control and Optimization} etc.
His current research interests include system modeling,
 adaptive
control, stochastic systems, and multi-agent systems.

Besides IEEE Fellow, Dr. Zhang is also an IFAC Fellow, CAA Fellow,
a Member of the European Academy of Sciences and Arts, and an
Academician of the International Academy for Systems and Cybernetic
Sciences. He received twice the Second Prize of the State Natural
Science Award of China in 2010 and 2015, respectively, the
 Distinguished
Young Scholar Fund from National Natural Science
Foundation of China
in 1997, the First Prize of the Young Scientist
 Award of CAS in 1995,
the Outstanding Advisor Award of CAS in 2007,
2008 and 2009,
respectively. He is a Vice-Chair of the
IFAC Technical
Board, and a Convenor of Systems Science Discipline, Academic
 Degree Committee of
the State Council of China; and was a member of the Board of Governors,
IEEE Control
Systems Society; Vice President of
 the Systems Engineering
Society of China, Vice President of
 the Chinese
 Association of Automation,
General Co-Chair of the 33rd
and the 36th Chinese Control Conferences;
IPC Chair of the 2012 IEEE
Conference on Control
Applications, the 9th World Congress on
 Intelligent Control
and Automation, and the 17th IFAC Symposium on
 System Identification;
and is an IPC Vice-Chair of the 20th IFAC
World Congress, etc.

\end{IEEEbiography}

\end{document}